\documentclass[smallextended,envcountsect,]{svjour3}

\smartqed

\usepackage{amsfonts,amssymb,amsmath,mathtools}

\newcommand{\AC}{\operatorname{AC}}
\newcommand{\Cont}{\operatorname{C}}
\newcommand{\rd}{\mathrm{d}}
\newcommand{\dist}{\operatorname{dist}}
\newcommand{\rpm}{\operatorname{rpm}}
\newcommand{\co}{\operatorname{co}}
\newcommand{\Id}{\operatorname{Id}_n}
\newcommand{\diam}{\operatorname{diam}}
\newcommand{\argmin}[1]{\underset{#1}{\operatorname{argmin}} \;}
\newcommand{\sgn}{\operatorname{sgn}}

\begin{document}

\title{Sensitivity Analysis of Value Functional of Fractional Optimal Control Problem with Application to Construction of Optimal Feedback Control}

\author{Mikhail Gomoyunov}

\institute{Mikhail Gomoyunov \at
    N.N. Krasovskii Institute of Mathematics and Mechanics of the Ural Branch of the Russian Academy of Sciences,
    Ekaterinburg 620108, Russia\\
    Ural Federal University,
    Ekaterinburg 620002, Russia \\
    m.i.gomoyunov@gmail.com
}

\date{Received: date / Accepted: date}

\maketitle

\begin{abstract}
    We consider an optimal control problem for a dynamical system described by a Caputo fractional differential equation and a terminal cost functional.
    We prove that, under certain assumptions, the (non-smooth, in general) value functional of this problem has a property of directional differentiability of order $\alpha$.
    As an application of this result, we propose a new method for constructing an optimal positional (feedback) control strategy.
\end{abstract}
\keywords{Optimal control problem \and Fractional differential equation \and Value functional \and Directional derivative \and Feedback control \and Hamilton--Jacobi--Bellman equation}
\subclass{49N35 \and 34A08 \and 49L12 \and 49J52}

\section{Introduction}
\label{section_Introduction}

    In this paper, we deal with a finite-horizon optimal control problem involving a dynamical system described by a non-linear differential equation with Caputo fractional derivative of order $\alpha \in (0, 1)$ and a terminal cost functional to be minimized.
    We focus on the question related to the construction of optimal positional control strategies, which allow us to generate $\varepsilon$-optimal controls for any predetermined accuracy $\varepsilon > 0$ by using the corresponding time-discrete recursive feedback control procedures (the terminology is borrowed from the theory of positional differential games, see, e.g., \cite{Krasovskii_Subbotin_1988,Krasovskii_Krasovskii_1995} and also, e.g., \cite{Osipov_1971,Lukoyanov_2000_JAMM,Lukoyanov_2003_2} for the case of time-delay systems).
    We follow the approach that goes back to the dynamic programming principle and relies on the consideration of the value functional (the functional of optimal result) and the analysis of its properties.
    Let us recall and briefly discuss the previous results obtained in this direction, which are the background of the present paper.

    In \cite{Gomoyunov_2020_SIAM} (see also \cite{Gomoyunov_2021_Mathematics}), it was proposed to associate the optimal control problem under consideration with a Cauchy problem for the corresponding Hamilton--Jacobi--Bellman equation with so-called fractional coinvariant derivatives ($ci$-derivatives for short) of order $\alpha$ and the natural right-end boundary condition (the use of the terminology of $ci$-derivatives goes back to \cite{Kim_1999,Lukoyanov_2000_JAMM}, see also the discussion in \cite[Section 5.2]{Gomoyunov_Lukoyanov_Plaksin_2021}).
    In particular, it was shown that, if the value functional is $ci$-smooth of order $\alpha$ (i.e., if it is continuous and has continuous $ci$-derivatives of order $\alpha$), then it can be characterized as a unique (classical) solution of the associated Cauchy problem and, moreover, an optimal positional control strategy can be constructed by extremal aiming in the direction of the $ci$-gradient of order $\alpha$ of this functional \cite[Theorems 10.1 and 11.1]{Gomoyunov_2020_SIAM}.
    However, by analogy with the case of optimal control problems for dynamical systems described by ordinary differential equations (i.e., when $\alpha = 1$), the value functional usually does not have the required smoothness properties.
    Furthermore, it is not clear how these properties of the value functional can be ensured by imposing assumptions on the input data of the problem (i.e., on the function from the right-hand side of the dynamic equation and on the function determining the terminal cost functional).

    On the other hand, in a fairly general non-smooth case, it was established in \cite[Theorem 1]{Gomoyunov_Lukoyanov_2021} (see also \cite[Theorem 2]{Gomoyunov_2021_Mathematics}) that the value functional of the optimal control problem coincides with a unique generalized (in the minimax sense) solution of the associated Cauchy problem (for details on the theory of minimax solutions of Hamilton--Jacobi equations with partial derivatives and with $ci$-derivatives of first order (i.e., when $\alpha = 1$), the reader is referred to, e.g., \cite{Subbotin_1995} and \cite{Gomoyunov_Lukoyanov_Plaksin_2021}, respectively).
    In particular, a pair of differential inequalities was obtained that characterizes the (non-smooth, in general) value functional in terms of so-called lower and upper derivatives of order $\alpha$ in (single-valued) directions \cite[Theorem 2 and formula (8.8)]{Gomoyunov_Lukoyanov_2021}.
    In addition, in \cite[Section 6]{Gomoyunov_2019_Trudy_Eng} and \cite[Section 9]{Gomoyunov_2021_Mathematics}, two general methods for constructing optimal positional control strategies based on the (non-smooth, in general) value functional were given.
    Nevertheless, let us emphasize that, even if the value functional has already been found, both of these methods require solving rather complex infinite-dimensional constrained optimization problems to select the desired extremal directions, which makes them difficult to apply even in relatively simple examples.
    For completeness, let us also mention two other possible approaches to constructing optimal positional control strategies: the approximation schemes from \cite{Gomoyunov_2020_DGA} and, in the case of linear dynamics, the technique from \cite{Gomoyunov_2020_ACS} based on reducing the optimal control problem to an auxiliary optimal control problem for a dynamical system described by an ordinary differential equation.

    In the present paper, we propose a new method for constructing an optimal positional control strategy, which devoid of the indicated disadvantages of the results of \cite{Gomoyunov_2019_Trudy_Eng,Gomoyunov_2020_SIAM,Gomoyunov_2021_Mathematics}.
    More precisely, in contrast to \cite{Gomoyunov_2020_SIAM}, we deal with the case when the value functional of the optimal control problem under consideration can be non-smooth and formulate all the necessary assumptions in terms of the input data of the problem.
    At the same time, compared to \cite{Gomoyunov_2019_Trudy_Eng,Gomoyunov_2021_Mathematics}, the new method is simpler and more efficient, since it requires only the calculation of directional derivatives of order $\alpha$ of the value functional.

    The proposed method can be considered as a modification of the method developed for the case of dynamical systems described by ordinary differential equations in, e.g., \cite[Theorem II.15]{Subbotina_2006}.
    In this connection, the reader is also referred to, e.g., \cite{Berkovitz_1989,Frankowska_1989,Rowland_Vinter_1991}.
    The basis of the method is the characterization of the value functional obtained in \cite[Theorem 2 and formula (8.8)]{Gomoyunov_Lukoyanov_2021} and a property of directional differentiability of order $\alpha$ of the value functional.
    The proof of this property, which is also of independent interest, is the main part of the present paper.
    Despite the fact that it is carried out mainly according to the scheme from, e.g., \cite{Subbotin_Subbotina_1982_DAN} and \cite[Theorem 1]{Subbotin_Subbotina_1983} (see also \cite[Theorem II.6]{Subbotina_2006}), it fully takes into account the features associated with the fact that the dynamical system is described by the Caputo fractional differential equation.
    Let us briefly outline the main steps of the proof.
    At the first step, we pass to the use of relaxed controls (also sometimes called generalized controls).
    Since the right-hand side of the dynamic equation has the same form as in the case of ordinary differential equations, we borrow the notion of a relaxed control and the corresponding technique from the optimal control theory for the ordinary case (see, e.g., \cite[Chapter IV]{Warga_1972}, \cite[Part 3]{Fattorini_1999}, and also \cite[Section 6.1]{Krasovskii_Subbotin_1988}).
    At the second step, we make a change of the time variable in order to obtain the unified control interval $[0, 1]$ for all possible initial times, and, thus, we move to some auxiliary (weakly-singular) Volterra integral equation.
    The third step is to analyze the dependence of a solution of the auxiliary integral equation on the parameters, which are the initial data and the relaxed control.
    In particular, we study directional differentiability of order $\alpha$ of the endpoint of the solution with respect to initial data.
    Note that these three steps are essentially based on the results of \cite{Gomoyunov_2022_FCAA}, which, however, need to be further developed since now we are dealing with the differential equation containing the control and some uniform estimates are required.
    At the final step, we prove a theorem on directional differentiability of order $\alpha$ of a lower envelope of a family of uniformly directionally differentiable of order $\alpha$ functionals and apply it to our setting.
    This theorem is an analog of, e.g., \cite[Proposition 3]{Lukoyanov_2001_PMM_Eng}, where a close result was obtained for the case of first-order $ci$-derivatives.

    The rest of the paper is organized as follows.
    In Section \ref{section_OCP}, we formulate the optimal control problem, which is the subject of the present study.
    After some preliminaries, we describe the dynamical system and the cost functional, introduce the value functional, and give a definition of an optimal positional control strategy.
    In Section \ref{section_HJB}, we consider the associated Cauchy problem for the Hamilton--Jacobi--Bellman equation with $ci$-derivatives of order $\alpha$ and the right-end boundary condition and briefly describe some of the results obtained in \cite{Gomoyunov_2020_SIAM,Gomoyunov_Lukoyanov_2021}.
    In particular, we recall the notions of lower and upper directional derivatives of order $\alpha$ of functionals.
    Sections \ref{section_Envelopes}--\ref{section_Directional_differentiability} are devoted to the detailed proof of the property of directional differentiability of order $\alpha$ of the value functional (see Theorem \ref{theorem_directional_differentiability}).
    In Section \ref{section_optimal_control}, we describe the new method for constructing an optimal positional control strategy (see Theorem \ref{theorem_optimal_control}).
    The paper concludes with an illustrative example, given in Section \ref{section_Example}.

\section{Optimal Control Problem}
\label{section_OCP}

    \subsection{Preliminaries}

        Let $\alpha \in (0, 1)$, $n \in \mathbb{N}$, and $T > 0$ be fixed throughout the paper.
        Let $\mathbb{R}^n$ be the Euclidean space of $n$-dimensional vectors with the inner product $\langle \cdot, \cdot \rangle$ and the norm $\|\cdot\|$, and let $\mathbb{R}^{n \times n}$ be the space of ($n \times n$)-dimensional matrices endowed with the corresponding induced norm, also denoted by $\|\cdot\|$.

        For every $t \in [0, T]$, let us introduce the space $\AC^\alpha([0, t], \mathbb{R}^n)$ of all functions $x \colon [0, t] \to \mathbb{R}^n$ each of which can be represented in the following form for some (Lebesgue) measurable and essentially bounded function $f \colon [0, t] \to \mathbb{R}^n$ (see, e.g., \cite[Definition 2.3]{Samko_Kilbas_Marichev_1993}):
        \begin{equation} \label{x_f}
            x(\tau)
            = x(0) + \frac{1}{\mathrm{\Gamma}(\alpha)} \int_{0}^{\tau} \frac{f(\xi)}{(\tau - \xi)^{1 - \alpha}} \, \rd \xi
            \quad \forall \tau \in [0, t].
        \end{equation}
        Note that, in the right-hand side of the equality from \eqref{x_f}, the second term is the {\it Riemann--Liouville fractional integral of order} $\alpha$ of the function $f(\cdot)$ (see, e.g., \cite[Definition 2.1]{Samko_Kilbas_Marichev_1993}) and $\mathrm{\Gamma}$ is the gamma-function.
        In accordance with, e.g., \cite[Remark 3.3]{Samko_Kilbas_Marichev_1993}, we consider the space $\AC^\alpha([0, t], \mathbb{R}^n)$ as a subset of the space $\Cont([0, t], \mathbb{R}^n)$ of all continuous functions from $[0, t]$ to $\mathbb{R}^n$ endowed with the uniform (supremum) norm:
        \begin{equation*}
            \|x(\cdot)\|_{[0, t]}
            \coloneqq \max_{\tau \in [0, t]} \|x(\tau)\|
            \quad \forall x(\cdot) \in \Cont([0, t], \mathbb{R}^n).
        \end{equation*}
        In addition, due to, e.g., \cite[Theorem 2.4]{Samko_Kilbas_Marichev_1993}, every function $x(\cdot) \in \AC^\alpha([0, t], \mathbb{R}^n)$ has at almost every (a.e.) $\tau \in [0, t]$ a {\it Caputo fractional derivative of order} $\alpha$, which is defined by (see, e.g., \cite[Section 2.4]{Kilbas_Srivastava_Trujillo_2006} and \cite[Chapter 3]{Diethelm_2010})
        \begin{equation} \label{Caputo}
            (^C D^\alpha x)(\tau)
            \coloneqq \frac{1}{\mathrm{\Gamma}(1 - \alpha)} \frac{\rd}{\rd \tau} \int_{0}^{\tau} \frac{x(\xi) - x(0)}{(\tau - \xi)^\alpha} \, \rd \xi.
        \end{equation}
        Moreover, if representation \eqref{x_f} is valid for some measurable and essentially bounded function $f(\cdot)$, then $(^C D^\alpha x)(\tau) = f(\tau)$ for a.e. $\tau \in [0, t]$.

    \subsection{Dynamic Equation}

        Let us consider a {\it dynamical system} described by the differential equation
        \begin{equation} \label{system}
            (^C D^\alpha x)(\tau)
            = f(\tau, x(\tau), u(\tau)).
        \end{equation}
        Here, $\tau \in [0, T]$ is time,
        $x(\tau) \in \mathbb{R}^n$ is the current value of the state vector,
        $(^C D^\alpha x)(\tau)$ is the Caputo fractional derivative of order $\alpha$ (see \eqref{Caputo}),
        $u(\tau) \in P$ is the current value of the control vector,
        and $P \subset \mathbb{R}^{n_u}$ is a compact set, $n_u \in \mathbb{N}$.
        Let us assume that the following conditions hold:
        \begin{description}
            \item[$(f.1)$]
                The function $f \colon [0, T] \times \mathbb{R}^n \times P \to \mathbb{R}^n$ is continuous.

            \item[$(f.2)$]
                For any $R \geq 0$, there exists $\lambda_f \geq 0$ such that, for any $\tau \in [0, T]$, any $x$, $x^\prime \in B(R)$, and any $u \in P$,
                \begin{equation*}
                    \|f(\tau, x, u) - f(\tau, x^\prime, u)\|
                    \leq \lambda_f \|x - x^\prime\|,
                \end{equation*}
                where we denote $B(R) \coloneqq \{ x \in \mathbb{R}^n \colon \|x\| \leq R\}$.

            \item[$(f.3)$]
                There exists $c_f \geq 0$ such that
                \begin{equation*}
                    \|f(\tau, x, u)\|
                    \leq c_f (1 + \|x\|)
                    \quad \forall \tau \in [0, T] \quad \forall x \in \mathbb{R}^n \quad \forall u \in P.
                \end{equation*}
        \end{description}

    \subsection{Space of Positions}

        In accordance with, e.g., \cite{Gomoyunov_2020_SIAM}, by a {\it position} of system \eqref{system}, we mean a pair $(t, w(\cdot))$ consisting of a time $t \in [0, T]$ and a function $w(\cdot) \in \AC^\alpha([0, t], \mathbb{R}^n)$, which is treated as a history of a motion of system \eqref{system} on the time interval $[0, t]$.
        Let us denote the set of all such positions $(t,w(\cdot ))$ by $G$, i.e.,
        \begin{equation*}
            G \coloneqq \bigcup_{t \in [0, T]} (\{t\} \times \AC^\alpha([0, t], \mathbb{R}^n)),
        \end{equation*}
        and also introduce the set $G^0 \coloneqq \{ (t, w(\cdot)) \in G \colon t < T \}$.
        Let us endow the set $G$ (and, respectively, its subset $G^0$) with the metric
        \begin{equation} \label{dist}
            \dist \bigl( (t, w(\cdot)), (t^\prime, w^\prime(\cdot)) \bigr)
            \coloneqq |t - t^\prime| + \max_{\tau \in [0, T]} \|w(\min\{\tau, t\}) - w^\prime(\min\{\tau, t^\prime\})\|
        \end{equation}
        for all $(t, w(\cdot))$, $(t^\prime, w^\prime(\cdot)) \in G$.
        Let us note that the mapping
        \begin{equation} \label{general_mapping}
            [0, T] \times \AC^\alpha([0, T], \mathbb{R}^n) \ni (t, x(\cdot)) \mapsto (t, x_t(\cdot)) \in G
        \end{equation}
        is continuous, where $x_t(\cdot) \in \AC^\alpha([0, t], \mathbb{R}^n)$ stands for the restriction of the function $x(\cdot)$ to the interval $[0, t]$, i.e.,
        \begin{equation} \label{x_t}
            x_t(\tau)
            \coloneqq x(\tau)
            \quad \forall \tau \in [0, t].
        \end{equation}

        \begin{remark}
            In \cite{Gomoyunov_2020_SIAM,Gomoyunov_2020_DE,Gomoyunov_Lukoyanov_2021}, the set of positions $G$ was considered with a different metric.
            However, despite the fact that this metric and the metric $\dist$ from \eqref{dist} are not strongly equivalent, they induce the same topology on $G$ (see, e.g., \cite[Section 5.1]{Gomoyunov_Lukoyanov_Plaksin_2021} for details).
            This allows us to apply the results from \cite{Gomoyunov_2020_SIAM,Gomoyunov_2020_DE,Gomoyunov_Lukoyanov_2021} in the present paper.
        \end{remark}

    \subsection{Open-Loop Controls and Motions of the System}
    \label{subsection_admissible_control_and_motions}

        Let an {\it initial position} $(t, w(\cdot)) \in G^0$ be given.
        Let us consider the set
        \begin{equation} \label{X}
            X(t, w(\cdot))
            \coloneqq \bigl\{ x(\cdot) \in \AC^\alpha([0, T], \mathbb{R}^n) \colon
            x_t(\cdot) = w(\cdot) \bigr\}.
        \end{equation}
        In addition, let us introduce the function
        \begin{equation} \label{a}
            a(\tau \mid t, w(\cdot))
            \coloneqq \begin{cases}
                w(\tau),
                & \mbox{if } \tau \in [0, t], \\
                \displaystyle
                w(0) + \frac{1}{\mathrm{\Gamma}(\alpha)} \int_{0}^{t} \frac{(^C D^\alpha w)(\xi)}{(\tau - \xi)^{1 - \alpha}} \, \rd \xi,
                & \mbox{if } \tau \in (t, T].
            \end{cases}
        \end{equation}
        Let us note that $a(\cdot \mid t, w(\cdot)) \in X(t, w(\cdot))$ and $(^C D^\alpha a(\cdot \mid t, w(\cdot)))(\tau) = 0$ for all $\tau \in (t, T)$.
        It is also convenient to formally put $a(\cdot \mid T, w(\cdot)) \coloneqq w(\cdot)$ for all $w(\cdot) \in \AC^\alpha([0, T], \mathbb{R}^n)$.

        By an {\it admissible} ({\it open-loop}) {\it control} on the time interval $[t, T]$, we mean any measurable function $u \colon [t, T] \to P$.
        Let $\mathcal{U}(t, T)$ be the set of all such controls.
        A {\it motion} of system \eqref{system} generated from the initial position $(t, w(\cdot))$ by a control $u(\cdot) \in \mathcal{U}(t, T)$ is defined as a function $x(\cdot) \in X(t, w(\cdot))$ that together with the function $u(\cdot)$ satisfies the differential equation \eqref{system} for a.e. $\tau \in [t, T]$.
        According to, e.g., \cite[Proposition 2]{Gomoyunov_2020_DGA}, such a motion $x(\cdot)$ exists and is unique, and we denote it by $x(\cdot) \coloneqq x(\cdot \mid t, w(\cdot), u(\cdot))$.
        Moreover, this motion $x(\cdot)$ is a unique function from $\Cont([0, T], \mathbb{R}^n)$ that satisfies the initial condition $x_t(\cdot) = w(\cdot)$ (see \eqref{x_t}) and the {\it Volterra integral equation}
        \begin{equation*}
            x(\tau)
            = a(\tau \mid t, w(\cdot))
            + \frac{1}{\mathrm{\Gamma}(\alpha)} \int_{t}^{\tau} \frac{f(\xi, x(\xi), u(\xi))}{(\tau - \xi)^{1 - \alpha}} \, \rd \xi
            \quad \forall \tau \in [t, T].
        \end{equation*}

    \subsection{Cost Functional}

        Let us consider the {\it terminal cost functional}
        \begin{equation} \label{cost_functional}
            J(t, w(\cdot), u(\cdot))
            \coloneqq \sigma\bigl( x(T \mid t, w(\cdot), u(\cdot)) \bigr)
            \ \ \forall (t, w(\cdot)) \in G^0 \ \ \forall u(\cdot) \in \mathcal{U}(t, T).
        \end{equation}
        For a given initial position $(t, w(\cdot)) \in G^0$, the problem is to find a control $u(\cdot) \in \mathcal{U}(t, T)$ that {\it minimizes} the value of this functional.
        Let us assume that the function $\sigma \colon \mathbb{R}^n \to \mathbb{R}$ satisfies the following condition:
        \begin{description}
            \item[$(\sigma.1)$]
                For any $R \geq 0$, there exists $\lambda_\sigma \geq 0$ such that
                \begin{equation*}
                    |\sigma(x) - \sigma(x^\prime)|
                    \leq \lambda_\sigma \|x - x^\prime\|
                    \quad \forall x, x^\prime \in B(R).
                \end{equation*}
        \end{description}

        Thus, we deal with the {\it optimal control problem} for the dynamical system \eqref{system} and the cost functional \eqref{cost_functional} under assumptions $(f.1)$--$(f.3)$ and $(\sigma.1)$.
        Let us note that the main results of the paper will be obtained under more restrictive assumptions on the functions $f$ and $\sigma$, which will be formulated when they are needed for the first time.

    \subsection{Value Functional and $\varepsilon$-Optimal Open-Loop Controls}

        Let us introduce the {\it value functional} $\rho \colon G \to \mathbb{R}$ of the optimal control problem \eqref{system}, \eqref{cost_functional} by
        \begin{align}
            \rho(t, w(\cdot))
            & \coloneqq \inf_{u(\cdot) \in \mathcal{U}(t, T)} J(t, w(\cdot), u(\cdot))
            \nonumber \\
            & = \inf_{u(\cdot) \in \mathcal{U}(t, T)} \sigma \bigl( x(T \mid t, w(\cdot), u(\cdot)) \bigr)
            \quad \forall (t, w(\cdot)) \in G^0
            \label{value_functional_G^0}
        \end{align}
        and, formally,
        \begin{equation} \label{value_functional_T}
            \rho(T, w(\cdot))
            \coloneqq \sigma(w(T))
            \quad \forall w(\cdot) \in \AC^\alpha([0, T], \mathbb{R}^n).
        \end{equation}
        Given an initial position $(t, w(\cdot)) \in G^0$ and a number $\varepsilon > 0$, we call a control $u(\cdot) \in \mathcal{U}(t, T)$ {\it $\varepsilon$-optimal} if
        \begin{equation*}
            J(t, w(\cdot), u(\cdot))
            = \sigma \bigl( x(T \mid t, w(\cdot), u(\cdot)) \bigr)
            \leq \rho(t, w(\cdot)) + \varepsilon.
        \end{equation*}

    \subsection{Positional Control Strategies}
    \label{section_positional}

        In accordance with, e.g., \cite[Section 11]{Gomoyunov_2020_SIAM}, by a {\it positional control strategy} in the optimal control problem \eqref{system}, \eqref{cost_functional}, we mean any mapping $U \colon G^0 \to P$.

        Let an initial position $(t, w(\cdot)) \in G^0$ be given, and let $\Delta$ be a {\it partition} of the time interval $[t, T]$, i.e.,
        \begin{equation*}
            \Delta
            \coloneqq \{\tau_j\}_{j \in \overline{1, k + 1}},
            \quad \tau_1 = t,
            \quad \tau_{j + 1}
            > \tau_j
            \quad \forall j \in \overline{1, k},
            \quad \tau_{k + 1} = T,
        \end{equation*}
        where $k \in \mathbb{N}$ and the notation $\overline{1, k + 1} \coloneqq \{j \in \mathbb{N} \colon j \leq k + 1\}$ is used.
        Then, a positional control strategy $U$ forms a piecewise constant control $u(\cdot) \in \mathcal{U}(t, T)$ and the corresponding motion $x(\cdot) \coloneqq x(\cdot \mid t, w(\cdot), u(\cdot))$ of system \eqref{system} by the following {\it recursive feedback control procedure}:
        \begin{equation*}
            u(\tau)
            \coloneqq U(\tau_j, x_{\tau_j}(\cdot))
            \quad \forall \tau \in [\tau_j, \tau_{j + 1}) \quad \forall j \in \overline{1, k},
        \end{equation*}
        where $x_{\tau_j}(\cdot)$ is the history of the motion $x(\cdot)$ on the time interval $[0, \tau_j]$ (see \eqref{x_t}).
        Formally putting $u(T) \coloneqq \tilde{u}$ for some fixed $\tilde{u} \in P$, we conclude that the described control procedure determines $u(\cdot)$ and $x(\cdot)$ uniquely.
        Let us denote the obtained (open-loop) control by $u(\cdot) \coloneqq u(\cdot \mid t, w(\cdot), U, \Delta)$.

        A positional control strategy $U$ is called {\it optimal} if the following holds:
        for any $(t, w(\cdot)) \in G^0$ and any $\varepsilon > 0$, there exists $\delta > 0$ such that, for any partition $\Delta \coloneqq \{\tau_j\}_{j \in \overline{1, k + 1}}$ with the diameter $\diam(\Delta) \coloneqq \max_{j \in \overline{1, k}} (\tau_{j + 1} - \tau_j) \leq \delta$, the control $u(\cdot \mid t, w(\cdot), U, \Delta)$ is $\varepsilon$-optimal, i.e.,
        \begin{equation*}
            \sigma \bigl( x(T \mid t, w(\cdot), u(\cdot \mid t, w(\cdot), U, \Delta)) \bigr)
            \leq \rho(t, w(\cdot)) + \varepsilon.
        \end{equation*}

\section{Hamilton--Jacobi--Bellman Equation}
\label{section_HJB}

    \subsection{Smooth Case}

        Following, e.g., \cite[Section 9]{Gomoyunov_2020_SIAM}, we say that a functional $\varphi \colon G \to \mathbb{R}$ is {\it $ci$-differ\-entiable of order} $\alpha$ at a point $(t, w(\cdot)) \in G^0$ if there are $\partial_t^\alpha \varphi(t, w(\cdot)) \in \mathbb{R}$ and $\nabla^\alpha \varphi(t, w(\cdot)) \in \mathbb{R}^n$ such that, for any function $x(\cdot) \in X(t, w(\cdot))$ (see \eqref{X}),
        \begin{align*}
            & \biggl| \frac{\varphi(t + \delta, x_{t + \delta}(\cdot)) - \varphi(t, w(\cdot))}{\delta} \\
            & \quad - \partial_t^\alpha \varphi(t, w(\cdot))
            - \biggl\langle \nabla^\alpha \varphi(t, w(\cdot)), \frac{1}{\delta} \int_{t}^{t + \delta} (^C D^\alpha x)(\xi) \, \rd \xi \biggr\rangle \biggr|
            \to 0
            \quad \text{as } \delta \to 0^+,
        \end{align*}
        where $x_{t + \delta}(\cdot)$ is the restriction of the function $x(\cdot)$ to the interval $[0, t + \delta]$ (see \eqref{x_t}) and $(^C D^\alpha x)(\xi)$ is the Caputo fractional derivative of order $\alpha$ of the function $x(\cdot)$ (see \eqref{Caputo}).
        In this case, the values $\partial_t^\alpha \varphi(t, w(\cdot))$ and $\nabla^\alpha \varphi(t, w(\cdot))$ are called respectively the {\it $ci$-derivative in $t$ of order} $\alpha$ and the {\it $ci$-gradient of order} $\alpha$ of the functional $\varphi$ at the point $(t, w(\cdot))$.
        Moreover, we say that a functional $\varphi \colon G \to \mathbb{R}$ is {\it $ci$-smooth of order} $\alpha$ if it is continuous, $ci$-differentiable of order $\alpha$ at every point $(t, w(\cdot)) \in G^0$, and the mappings $\partial_t^\alpha \varphi \colon G^0 \to \mathbb{R}$ and $\nabla^\alpha \varphi \colon G^0 \to \mathbb{R}^n$ are continuous.

        With the optimal control problem \eqref{system}, \eqref{cost_functional}, let us associate the {\it Hamiltonian}
        \begin{equation} \label{Hamiltonian}
            H(\tau, x, s)
            \coloneqq \min_{u \in P} \langle s, f(\tau, x, u) \rangle
            \quad \forall \tau \in [0, T] \quad \forall x, s \in \mathbb{R}^n
        \end{equation}
        and the {\it Cauchy problem} for the {\it Hamilton--Jacobi--Bellman equation}
        \begin{equation} \label{HJB}
            \partial_t^\alpha \varphi(t, w(\cdot)) + H\bigl(t, w(t), \nabla^\alpha \varphi(t, w(\cdot)) \bigr)
            = 0
            \quad \forall (t, w(\cdot)) \in G^0
        \end{equation}
        under the right-end {\it boundary condition}
        \begin{equation} \label{boundary_condition}
            \varphi(T, w(\cdot))
            = \sigma(w(T))
            \quad \forall w(\cdot) \in \AC^\alpha([0, T], \mathbb{R}^n),
        \end{equation}
        where a functional $\varphi \colon G \to \mathbb{R}$ is the unknown.

        By \cite[Theorems 10.1 and 11.1]{Gomoyunov_2020_SIAM}, we have the following two results.
        The first one is a criteria for a $ci$-smooth of order $\alpha$ functional to be the value functional.
        The second one proposes a method for constructing an optimal positional control strategy in such a smooth case.
        \begin{theorem} \label{theorem_smooth_1}
            Under assumptions $(f.1)$--$(f.3)$ and $(\sigma.1)$, a $ci$-smooth of order $\alpha$ functional $\varphi \colon G \to \mathbb{R}$ is the value functional of the optimal control problem \eqref{system}, \eqref{cost_functional} if and only if $\varphi$ satisfies the Hamilton--Jacobi--Bellman equation \eqref{HJB} and the boundary condition \eqref{boundary_condition}.
        \end{theorem}
        \begin{theorem} \label{theorem_smooth_2}
            Suppose that assumptions $(f.1)$--$(f.3)$ and $(\sigma.1)$ hold and that the value functional $\rho$ of the optimal control problem \eqref{system}, \eqref{cost_functional} is $ci$-smooth of order $\alpha$.
            Then, a positional control strategy $U^\circ$ satisfying the condition
            \begin{equation} \label{U^0_smooth}
                U^\circ(t, w(\cdot))
                \in \argmin{u \in P} \langle \nabla^\alpha \rho(t, w(\cdot)), f(t, w(t), u) \rangle
                \quad \forall (t, w(\cdot)) \in G^0
            \end{equation}
            is optimal in this problem.
        \end{theorem}

        However, the value functional $\rho$ is usually not $ci$-smooth of order $\alpha$, which narrows the applicability of Theorems \ref{theorem_smooth_1} and \ref{theorem_smooth_2}.
        Moreover, the $ci$-gradient of order $\alpha$ of the value functional $\rho$ may fail to exist at some points $(t, w(\cdot)) \in G^0$, and, therefore, relation \eqref{U^0_smooth} cannot be directly used to construct an optimal positional control strategy.
        The method for constructing an optimal positional control strategy proposed in the present paper can be considered as a modification of relation \eqref{U^0_smooth} that allows us to handle such situations.

    \subsection{Non-Smooth Case: Inequalities for Directional Derivatives of Order $\alpha$}

        For any $(t, w(\cdot)) \in G^0$ and any $f \in \mathbb{R}^n$, let us introduce the function
        \begin{equation} \label{x^f}
            x^{(f)}(\tau \mid t, w(\cdot))
            \coloneqq
            \begin{cases}
                w(\tau), & \mbox{if } \tau \in [0, t], \\
                \displaystyle
                a(\tau \mid t, w(\cdot)) + \frac{(\tau - t)^\alpha f}{\mathrm{\Gamma}(\alpha + 1)}, & \mbox{if } \tau \in (t, T],
            \end{cases}
        \end{equation}
        where the function $a(\cdot \mid t, w(\cdot))$ is given by \eqref{a}.
        Let us note that the inclusion $x^{(f)}(\cdot \mid t, w(\cdot)) \in X(t, w(\cdot))$ holds (see \eqref{X}) and, for any $\tau \in (t, T)$, the equality $(^C D^\alpha x^{(f)}(\cdot \mid t, w(\cdot)))(\tau) = f$ is valid.

        According to, e.g., \cite[Section 7]{Gomoyunov_Lukoyanov_2021}, the {\it lower} and {\it upper} ({\it right}) {\it derivatives of order} $\alpha$ of a functional $\varphi \colon G \to \mathbb{R}$ at a point $(t, w(\cdot)) \in G^0$ in a {\it direction} $f \in \mathbb{R}^n$ are defined as follows:
        \begin{align}
            \partial^\alpha_- \{\varphi(t, w(\cdot)) \mid f\}
            & \coloneqq \liminf_{\delta \to 0^+}
            \frac{\varphi(t + \delta, x^{(f)}_{t + \delta}(\cdot \mid t, w(\cdot))) - \varphi(t, w(\cdot))}{\delta},
            \nonumber \\
            \partial^\alpha_+ \{\varphi(t, w(\cdot)) \mid f\}
            & \coloneqq \limsup_{\delta \to 0^+}
            \frac{\varphi(t + \delta, x^{(f)}_{t + \delta}(\cdot \mid t, w(\cdot))) - \varphi(t, w(\cdot))}{\delta},
            \label{lower_and_upper_derivatives}
        \end{align}
        where $x^{(f)}_{t + \delta}(\cdot \mid t, w(\cdot))$ is the restriction of $x^{(f)}(\cdot \mid t, w(\cdot))$ to $[0, t + \delta]$ (see \eqref{x_t}).

        For a functional $\varphi \colon G \to \mathbb{R}$, let us consider the following {\it local Lipschitz continuity property}:
        \begin{description}
            \item[$(L)$]
                For every compact set $K \subset G$, there exists $\lambda_\varphi \geq 0$ such that
                \begin{equation*}
                    |\varphi(t, w(\cdot)) - \varphi(t, w^\prime(\cdot))|
                    \leq \lambda_\varphi \biggl(\|a(T) - a^\prime(T)\|
                    + \int_{t}^{T} \frac{\|a(\xi) - a^\prime(\xi)\|}{(T - \xi)^{1 - \alpha}} \, \rd \xi \biggr)
                \end{equation*}
                for all $(t, w(\cdot))$, $(t, w^\prime(\cdot)) \in K$, where the functions $a(\cdot) \coloneqq a(\cdot \mid t, w(\cdot))$ and $a^\prime(\cdot) \coloneqq a(\cdot \mid t, w^\prime(\cdot))$ are given by \eqref{a}.
        \end{description}

        Formally embedding the optimal control problem \eqref{system}, \eqref{cost_functional} into the corresponding zero-sum differential game with a fictitious second player (see, e.g., \cite[Section 5.1]{Subbotina_2006} for details) and then applying \cite[Theorem 2]{Gomoyunov_Lukoyanov_2021}, which should be slightly modified in view of the form of the cost functional \eqref{cost_functional}, we obtain (see also \cite[formula (8.8)]{Gomoyunov_Lukoyanov_2021})
        \begin{theorem} \label{theorem_non-smooth_general}
            Let assumptions $(f.1)$--$(f.3)$ and $(\sigma.1)$ hold.
            Then, a functional $\varphi \colon G \to \mathbb{R}$ is the value functional of the optimal control problem \eqref{system}, \eqref{cost_functional} if and only if $\varphi$ is continuous, possesses property $(L)$, meets the boundary condition \eqref{boundary_condition}, and satisfies the pair of differential inequalities
            \begin{equation*}
                \min_{f \in \co f(t, w(t), P)} \partial_-^\alpha \{\varphi(t, w(\cdot)) \mid f\}
                \leq 0
                \leq \min_{u \in P} \partial_+^\alpha \{\varphi(t, w(\cdot)) \mid f(t, w(t), u)\}
            \end{equation*}
            for all $(t, w(\cdot)) \in G^0$, where $\co f(t, w(t), P)$ denotes the convex hull of the set
            \begin{equation} \label{f(P)}
                f(t, w(t), P)
                \coloneqq \bigl\{ f(t, w(t), u) \in \mathbb{R}^n \colon
                u \in P \bigr\}.
            \end{equation}
        \end{theorem}

        For constructing an optimal positional control strategy in the problem \eqref{system}, \eqref{cost_functional} under assumptions $(f.1)$--$(f.3)$ and $(\sigma.1)$, the methods from \cite[Section 6]{Gomoyunov_2019_Trudy_Eng} and \cite[Section 9]{Gomoyunov_2021_Mathematics} can be applied.
        Each of these methods gives a formula for the optimal positional control strategy that has the same form as \eqref{U^0_smooth} but with $\nabla^\alpha \rho(t, w(\cdot))$ replaced by some {\it extremal direction} $s^\circ(t, w(\cdot))$.
        To determine the direction $s^\circ(t, w(\cdot))$, it is necessary to find a minimum point of the value functional $\rho$ as in \cite[Section 6]{Gomoyunov_2019_Trudy_Eng} or of the value functional $\rho$ perturbed by a certain Lyapunov--Krasovskii functional as in \cite[Section 9]{Gomoyunov_2021_Mathematics} over some specific compact set $\Theta(t, w(\cdot)) \subset \AC^\alpha([0, t], \mathbb{R}^n)$.
        The method  for constructing an optimal positional control strategy proposed in the present paper does not require solving any optimization problems to select the desired extremal directions.
        The basis of the method is Theorem \ref{theorem_non-smooth_general} and a property of directional differentiability of order $\alpha$ of the value functional $\rho$.

        We say that a functional $\varphi \colon G \to \mathbb{R}$ is {\it differentiable of order} $\alpha$ at a point $(t, w(\cdot)) \in G^0$ in a {\it direction} $f \in \mathbb{R}^n$ if the lower $\partial^\alpha_- \{\varphi(t, w(\cdot)) \mid f\}$ and upper $\partial^\alpha_+ \{\varphi(t, w(\cdot)) \mid f\}$ directional derivatives of order $\alpha$ are finite and coincide with each other.
        In this case, the {\it derivative of order} $\alpha$ of the functional $\varphi$ at the point $(t, w(\cdot))$ in the {\it direction} $f$ is given by
        \begin{equation*}
            \partial^\alpha \{\varphi(t, w(\cdot)) \mid f\}
            \coloneqq \lim_{\delta \to 0^+} \frac{\varphi(t + \delta, x^{(f)}_{t + \delta}(\cdot \mid t, w(\cdot))) - \varphi(t, w(\cdot))}{\delta}.
        \end{equation*}
        Let us note that, if the functional $\varphi$ is $ci$-differentiable of order $\alpha$ at the point $(t, w(\cdot))$, then it is differentiable of order $\alpha$ at this point in every direction $f \in \mathbb{R}^n$ and the equality below holds:
        \begin{equation} \label{directional_derivative_via_ci-derivatives}
            \partial^\alpha \{\varphi(t, w(\cdot)) \mid f\}
            = \partial_t^\alpha \varphi(t, w(\cdot)) + \langle \nabla^\alpha \varphi(t, w(\cdot)), f \rangle.
        \end{equation}
        In addition, we say that a functional $\varphi \colon G \to \mathbb{R}$ is {\it directionally differentiable of order} $\alpha$ if it is differentiable of order $\alpha$ at every point $(t, w(\cdot)) \in G^0$ and in every direction $f \in \mathbb{R}^n$.

        Thus, our first goal is to prove that, under some additional assumptions on the functions $f$ and $\sigma$, the value functional $\rho$ of the optimal control problem \eqref{system}, \eqref{cost_functional} is directionally differentiable of order $\alpha$.
        The scheme of the proof is borrowed from, e.g., \cite{Subbotin_Subbotina_1982_DAN} and \cite[Theorem 1]{Subbotin_Subbotina_1983} (see also \cite[Theorem II.6]{Subbotina_2006}) and is briefly outlined in Section \ref{section_Introduction} above.
        It seems convenient to begin with a general result on directional differentiability of order $\alpha$ of a lower envelope of a family of uniformly directionally differentiable of order $\alpha$ functionals, which will be applied at the final step of the proof.

\section{Directional Differentiability of Order $\alpha$ of Envelope of Family of Functionals}
\label{section_Envelopes}

    Let us consider a compact metric space $\mathcal{L}$ (the {\it space of parameters}) and a continuous functional $\psi \colon G^0 \times \mathcal{L} \to \mathbb{R}$.
    Let us assume that the functional $\psi$ possesses the following {\it property of uniform directional differentiability of order} $\alpha$: for any $(t, w(\cdot)) \in G^0$ and any $\ell \in \mathcal{L}$, there exist $\partial^\alpha_t \psi(t, w(\cdot), \ell) \in \mathbb{R}$ and $\nabla^\alpha \psi(t, w(\cdot), \ell) \in \mathbb{R}^n$ such that, for every direction $f \in \mathbb{R}^n$,
    \begin{align}
        & \biggl| \frac{\psi(t + \delta, x_{t + \delta}^{(f)}(\cdot \mid t, w(\cdot)), \ell) - \psi(t, w(\cdot), \ell)}{\delta}
        \nonumber \\
        & \quad - \partial^\alpha_t \psi(t, w(\cdot), \ell) - \langle \nabla^\alpha \psi(t, w(\cdot), \ell), f \rangle \biggr|
        \to 0
        \quad \text{as } \delta \to 0^+
        \label{psi_uniform_coinvariant_differentiability}
    \end{align}
    uniformly in the parameter $\ell \in \mathcal{L}$, where $x^{(f)}_{t + \delta}(\cdot \mid t, w(\cdot))$ is the restriction of the function $x^{(f)}(\cdot \mid t, w(\cdot))$ (see \eqref{x^f}) to the interval $[0, t + \delta]$ (see \eqref{x_t}).
    Let us assume also that, for every $(t, w(\cdot)) \in G^0$, the mappings below are continuous:
    \begin{equation} \label{psi_coinvariant_derivatives}
        \mathcal{L} \ni \ell \mapsto \partial^\alpha_t \psi(t, w(\cdot), \ell) \in \mathbb{R},
        \quad \mathcal{L} \ni \ell \mapsto \nabla^\alpha \psi(t, w(\cdot), \ell) \in \mathbb{R}^n.
    \end{equation}

    Now, let us suppose that a functional $\varphi \colon G \to \mathbb{R}$ is given such that
    \begin{equation} \label{varphi_is_envelope}
        \varphi(t, w(\cdot))
        = \min_{\ell \in \mathcal{L}} \psi(t, w(\cdot), \ell)
        \quad \forall (t, w(\cdot)) \in G^0,
    \end{equation}
    i.e., the functional $\varphi$ (more precisely, its restriction to $G^0$) is the lower envelope of the family of functionals $G^0 \ni (t, w(\cdot)) \mapsto \psi(t, w(\cdot), \ell) \in \mathbb{R}$ parameterized by $\ell \in \mathcal{L}$.
    Let us note that, due to compactness of the space $\mathcal{L}$ and continuity of the functional $\psi$, for every $(t, w(\cdot)) \in G^0$, the set
    \begin{equation*}
        \mathcal{L}^\circ(t, w(\cdot))
        \coloneqq \bigl\{ \ell \in \mathcal{L} \colon
        \psi(t, w(\cdot), \ell) = \varphi(t, w(\cdot)) \bigr\}
    \end{equation*}
    is non-empty and compact.
    Moreover, the multivalued mapping
    \begin{equation} \label{P^0}
        G^0 \ni (t, w(\cdot)) \mapsto \mathcal{L}^\circ(t, w(\cdot)) \subset \mathcal{L}
    \end{equation}
    is upper semicontinuous, i.e., if $(t_i, w_i(\cdot)) \in G^0$ and $\ell_i \in \mathcal{L}^\circ(t_i, w_i(\cdot))$ for all $i \in \mathbb{N}$ are such that $(t_i, w_i(\cdot)) \to (t_0, w_0(\cdot)) \in G^0$ and $\ell_i \to \ell_0 \in \mathcal{L}$ as $i \to \infty$, then $\ell_0 \in \mathcal{L}^\circ(t_0, w_0(\cdot))$.
    \begin{theorem} \label{theorem_envelope_theorem}
        Under the assumptions made above in this section, the functional $\varphi$ is directionally differentiable of order $\alpha$.
        In addition, for every point $(t, w(\cdot)) \in G^0$ and every direction $f \in \mathbb{R}^n$, the following equality is valid:
        \begin{equation} \label{directional_derivative_of_envelope}
            \partial^\alpha \{ \varphi(t, w(\cdot)) \mid f \}
            = \min_{\ell \in \mathcal{L}^\circ (t, w(\cdot))} \bigl( \partial^\alpha_t \psi(t, w(\cdot), \ell)
            + \langle \nabla^\alpha \psi(t, w(\cdot), \ell), f \rangle \bigr).
        \end{equation}
    \end{theorem}
    {\it Proof}
        We follow the scheme of the proof of \cite[Proposition 3]{Lukoyanov_2001_PMM_Eng}.

        Let $(t, w(\cdot)) \in G^0$ and $f \in \mathbb{R}^n$ be fixed.
        Let us note that the minimum in \eqref{directional_derivative_of_envelope} is attained owing to compactness of the set $\mathcal{L}^\circ (t, w(\cdot))$ and continuity of mappings \eqref{psi_coinvariant_derivatives}.

        In accordance with definition \eqref{lower_and_upper_derivatives} of the lower derivative of order $\alpha$ of the functional $\varphi$ at the point $(t, w(\cdot))$ in the direction $f$, let us choose a sequence $\{\delta_i\}_{i \in \mathbb{N}} \subset (0, T - t)$ such that $\delta_i \to 0$ and
        \begin{equation*}
            \frac{\varphi(t + \delta_i, w_i(\cdot)) - \varphi(t, w(\cdot))}{\delta_i}
            \to \partial^\alpha_- \{ \varphi(t, w(\cdot)) \mid f \}
            \quad \text{as } i \to \infty,
        \end{equation*}
        where we denote $w_i(\cdot) \coloneqq x^{(f)}_{t + \delta_i}(\cdot \mid t, w(\cdot))$ for all $i \in \mathbb{N}$.
        For every $i \in \mathbb{N}$, let us take an arbitrary $\ell_i \in \mathcal{L}^\circ(t + \delta_i, w_i(\cdot))$.
        Since the metric space $\mathcal{L}$ is compact, we can assume that there exists $\ell_0 \in \mathcal{L}$ such that $\ell_i \to \ell_0$ as $i \to \infty$.
        In addition, continuity of mapping \eqref{general_mapping} and the convergence $\delta_i \to 0$ as $i \to \infty$ imply that $(t + \delta_i, w_i(\cdot)) \to (t, w(\cdot))$ as $i \to \infty$.
        Hence, due to upper semicontinuity of the multivalued mapping \eqref{P^0}, we obtain that $\ell_0 \in \mathcal{L}^\circ(t, w(\cdot))$.
        Further, in view of \eqref{varphi_is_envelope}, for any $i \in \mathbb{N}$, we derive
        \begin{align*}
            & \frac{\varphi(t + \delta_i, w_i(\cdot)) - \varphi(t, w(\cdot))}{\delta_i}
            \geq \frac{\psi(t + \delta_i, w_i(\cdot), \ell_i) - \psi(t, w(\cdot), \ell_i)}{\delta_i} \\
            & \geq \partial^\alpha_t \psi(t, w(\cdot), \ell_i) + \langle \nabla^\alpha \psi(t, w(\cdot), \ell_i), f \rangle \\
            & - \biggl| \frac{\psi(t + \delta_i, w_i(\cdot), \ell_i) - \psi(t, w(\cdot), \ell_i)}{\delta_i}
            - \partial^\alpha_t \psi(t, w(\cdot), \ell_i) - \langle \nabla^\alpha \psi(t, w(\cdot), \ell_i), f \rangle \biggr|.
        \end{align*}
        Passing to the limit as $i \to \infty$, owing to continuity of derivatives \eqref{psi_coinvariant_derivatives} and the property of uniform directional differentiability of order $\alpha$ (see \eqref{psi_uniform_coinvariant_differentiability}), we get
        \begin{align}
            \partial^\alpha_- \{ \varphi(t, x(\cdot)) \mid f \}
            & \geq \partial^\alpha_t \psi(t, w(\cdot), \ell_0) + \langle \nabla^\alpha \psi(t, w(\cdot), \ell_0), f \rangle
            \nonumber \\
            & \geq \min_{\ell \in \mathcal{L}^\circ(t, w(\cdot))}
            \bigl( \partial^\alpha_t \psi(t, w(\cdot), \ell) + \langle \nabla^\alpha \psi(t, w(\cdot), \ell), f \rangle \bigr).
            \label{lower_derivative_geq}
        \end{align}

        On the other hand, using definition \eqref{lower_and_upper_derivatives} of the upper derivative of order $\alpha$ of $\varphi$ at $(t, w(\cdot))$ in the direction $f$, let us choose a sequence $\{\delta_i\}_{i \in \mathbb{N}} \subset (0, T - t)$ such that $\delta_i \to 0$ and
        \begin{equation*}
            \frac{\varphi(t + \delta_i, w_i(\cdot)) - \varphi(t, w(\cdot))}{\delta_i}
            \to \partial^\alpha_+ \{ \varphi(t, w(\cdot)) \mid f \}
            \quad \text{as } i \to \infty,
        \end{equation*}
        where we denote again $w_i(\cdot) \coloneqq x^{(f)}_{t + \delta_i}(\cdot \mid t, w(\cdot))$ for all $i \in \mathbb{N}$.
        Then, based on \eqref{psi_uniform_coinvariant_differentiability} and \eqref{varphi_is_envelope}, we derive
        \begin{align*}
            \partial^\alpha_+ \{ \varphi(t, w(\cdot)) \mid f \}
            & \leq \lim_{i \to \infty} \frac{\psi(t + \delta_i, w_i(\cdot), \ell) - \psi(t, w(\cdot), \ell)}{\delta_i} \\
            & = \partial^\alpha_t \psi(t, w(\cdot), \ell) + \langle \nabla^\alpha \psi(t, w(\cdot), \ell), f \rangle
        \end{align*}
        for all $\ell \in \mathcal{L}^\circ(t, w(\cdot))$.
        Hence, we have
        \begin{equation} \label{upper_derivative_leq}
            \partial^\alpha_+ \{ \varphi(t, w(\cdot)) \mid f \}
            \leq \min_{\ell \in \mathcal{L}^\circ(t, w(\cdot))}
            \bigl( \partial^\alpha_t \psi(t, w(\cdot), \ell) + \langle \nabla^\alpha \psi(t, w(\cdot), \ell), f \rangle \bigr).
        \end{equation}

        Since $\partial^\alpha_+ \{ \varphi(t, w(\cdot)) \mid f \} \geq \partial^\alpha_- \{ \varphi(t, w(\cdot)) \mid f \}$ due to \eqref{lower_and_upper_derivatives}, inequalities \eqref{lower_derivative_geq} and \eqref{upper_derivative_leq} imply that $\partial^\alpha_+ \{ \varphi(t, w(\cdot)) \mid f \} = \partial^\alpha_- \{ \varphi(t, w(\cdot)) \mid f \}$ and that formula \eqref{directional_derivative_of_envelope} takes place.
        The theorem is proved.
    \qed

    In order to apply Theorem \ref{theorem_envelope_theorem} to our setting, we need to represent the value functional $\rho$ of the optimal control problem \eqref{system}, \eqref{cost_functional} in form \eqref{varphi_is_envelope} for some space of parameters $\mathcal{L}$ and functional $\psi$ with the required properties.
    The desired representation is based directly on definition \eqref{value_functional_G^0} of $\rho$.
    Accordingly, one of the difficulties that arise here is related to the appropriate compactness property of the set of controls $\mathcal{U}(t, T)$ and the corresponding continuous dependence property of a motion $x(\cdot \mid t, w(\cdot), u(\cdot)) \in \AC^\alpha([0, T], \mathbb{R}^n)$ of system \eqref{system} with respect to a control $u(\cdot) \in \mathcal{U}(t, T)$.
    We overcome this difficulty by passing to the use of relaxed controls.

\section{Relaxed Controls}
\label{section_relaxed_controls}

    Since the right-hand side of the dynamic equation \eqref{system} is the same as in the case of ordinary differential equations (i.e., when $\alpha = 1$), we borrow the notion of a relaxed control from the optimal control theory for the ordinary case without any difference.
    The details can be found in, e.g., \cite[Chapter IV]{Warga_1972} and \cite[Part 3]{Fattorini_1999} (see also \cite[Section 6.1]{Krasovskii_Subbotin_1988}).

    Let $\rpm(P)$ be the set of all regular probability Borel measures on $P$ (recall that $P$ is a compact subset of $\mathbb{R}^{n_u}$ that describes the geometric constraints on the control in system \eqref{system}).
    The set $\rpm(P)$ is endowed with a metric such that, given a sequence of measures $\{\mu_i\}_{i \in \mathbb{N}} \subset \rpm(P)$ and a measure $\mu_0 \in \rpm(P)$, the convergence $\mu_i \to \mu_0$ as $i \to \infty$ means that
    \begin{equation} \label{convergence_rpm}
        \int_{P} \mathfrak{a}(u) \mu_i (\rd u)
        \to \int_{P} \mathfrak{a}(u) \mu_0 (\rd u)
        \quad \text{as } i \to \infty
    \end{equation}
    for all continuous functions $\mathfrak{a} \colon P \to \mathbb{R}$.
    The metric space $\rpm(P)$ is compact.

    Let us fix $t \in [0, T)$.
    A function $\mathfrak{b} \colon [t, T] \times P \to \mathbb{R}$ is called a {\it Carath\'{e}odory function} if it satisfies the following conditions:
    (i) for every $\tau \in [t, T]$, the function $P \ni u \mapsto \mathfrak{b}(\tau, u) \in \mathbb{R}$ is continuous;
    (ii) for every $u \in P$, the function $[t, T] \ni \tau \mapsto \mathfrak{b}(\tau, u) \in \mathbb{R}$ is measurable;
    (iii) there is an integrable function $\theta \colon [t, T] \to \mathbb{R}$ such that $|\mathfrak{b}(\tau, u)| \leq \theta(\tau)$ for all $\tau \in [t, T]$ and all $u \in P$.
    Let us denote by $\mathcal{M}(t, T)$ the set of all (equivalence classes of) measurable functions $\mu \colon [t, T] \to \rpm(P)$.
    Let us note that, for any Carath\'{e}odory function $\mathfrak{b} \colon [t, T] \times P \to \mathbb{R}$ and any function $\mu(\cdot) \in \mathcal{M}(t, T)$, the function
    \begin{equation*}
        \mathfrak{c}(\tau)
        \coloneqq \int_{P} \mathfrak{b}(\tau, u) \mu(\tau, \rd u)
        \quad \text{for a.e. } \tau \in [t, T]
    \end{equation*}
    is integrable.
    Here and below, we write $\mu(\tau, \rd u)$ instead of $\mu(\tau)(\rd u)$ for convenience.
    The set $\mathcal{M}(t, T)$ is endowed with a metric such that, given a sequence of functions $\{\mu_i(\cdot)\}_{i \in \mathbb{N}} \subset \mathcal{M}(t, T)$ and a function $\mu_0(\cdot) \in \mathcal{M}(t, T)$, the convergence $\mu_i(\cdot) \to \mu_0(\cdot)$ as $i \to \infty$ means that
    \begin{equation} \label{convergence_M}
        \int_{t}^{T} \int_{P} \mathfrak{b}(\tau, u) \mu_i(\tau, \rd u) \, \rd \tau
        \to \int_{t}^{T} \int_{P} \mathfrak{b}(\tau, u) \mu_0(\tau, \rd u) \, \rd \tau
        \quad \text{as } i \to \infty.
    \end{equation}
    for all Carath\'{e}odory functions $\mathfrak{b} \colon [t, T] \times P \to \mathbb{R}$.
    The metric space $\mathcal{M}(t, T)$ is compact.

    Any function $\mu(\cdot) \in \mathcal{M}(t, T)$ is considered as a {\it relaxed control} on the time interval $[t, T]$.
    In order to embed the set $\mathcal{U}(t, T)$ of (usual) controls into the space $\mathcal{M}(t, T)$, every function $u(\cdot) \in \mathcal{U}(t, T)$ is associated with the function $\mu_{u(\cdot)}(\cdot) \in \mathcal{M}(t, T)$ such that
    \begin{equation} \label{relaxed_control_via_usual_control}
        \mu_{u(\cdot)}(\tau)
        \coloneqq \delta_{u(\tau)}
        \quad \text{for a.e. } \tau \in [t, T],
    \end{equation}
    where $\delta_{u(\tau)} \in \rpm(P)$ is the Dirac measure at the point $u(\tau) \in P$.
    The set
    \begin{equation} \label{M_0}
        \mathcal{M}_\ast(t, T)
        \coloneqq \bigl\{ \mu_{u(\cdot)}(\cdot) \in \mathcal{M}(t, T) \colon
        u(\cdot) \in \mathcal{U}(t, T) \bigr\}
    \end{equation}
    is dense in $\mathcal{M}(t, T)$.

    By the function $f$ from the right-hand side of the dynamic equation \eqref{system}, let us define the function $f^\ast \colon [0, T] \times \mathbb{R}^n \times \rpm(P) \to \mathbb{R}^n$ as follows:
    \begin{equation} \label{f_ast}
        f^\ast(\tau, x, \mu)
        \coloneqq \int_{P} f(\tau, x, u) \mu(\rd u)
        \quad \forall \tau \in [0, T] \quad \forall x \in \mathbb{R}^n \quad \forall \mu \in \rpm(P).
    \end{equation}
    Then, due to assumptions $(f.1)$--$(f.3)$ on the function $f$ and definition \eqref{convergence_rpm} of convergence in the space $\rpm(P)$, the function $f^\ast$ is well-defined and has the properties listed below:
    \begin{description}
        \item[$(f^\ast.1)$]
            The function $f^\ast$ is continuous.

        \item[$(f^\ast.2)$]
            For any $R \geq 0$, there exists $\lambda_f \geq 0$ such that
            \begin{equation*}
                \|f^\ast(\tau, x, \mu) - f^\ast(\tau, x^\prime, \mu)\|
                \leq \lambda_f \|x - x^\prime\|
            \end{equation*}
            for all $\tau \in [0, T]$, all $x$, $x^\prime \in B(R)$, and all $\mu \in \rpm(P)$.

        \item[$(f^\ast.3)$]
            There exists $c_f \geq 0$ such that
            \begin{equation*}
                \|f^\ast(\tau, x, \mu)\|
                \leq c_f (1 + \|x\|)
                \quad \forall \tau \in [0, T] \quad \forall x \in \mathbb{R}^n \quad \forall \mu \in \rpm(P).
            \end{equation*}
    \end{description}
    Let us note that we keep the notation for the numbers $\lambda_f$ and $c_f$ in properties $(f^\ast.2)$ and $(f^\ast.3)$ since these numbers can be taken the same as in conditions $(f.2)$ and $(f.3)$.

    A {\it motion} of system \eqref{system} generated from an initial position $(t, w(\cdot)) \in G^0$ by a relaxed control $\mu(\cdot) \in \mathcal{M}(t, T)$ is defined as a function $x(\cdot) \in X(t, w(\cdot))$ (see \eqref{X}) such that $(^C D^\alpha x)(\tau) = f^\ast(\tau, x(\tau), \mu(\tau))$ for a.e. $\tau \in [t, T]$.
    As in the case of (usual) controls $u(\cdot) \in \mathcal{U}(t, T)$ (see Section \ref{subsection_admissible_control_and_motions}), it can be verified that such a motion $x(\cdot) \coloneqq x(\cdot \mid t, w(\cdot), \mu(\cdot))$ exists and is unique.
    Moreover, this motion $x(\cdot)$ is a unique function from $\Cont([0, T], \mathbb{R}^n)$ that satisfies the initial condition $x_t(\cdot) = w(\cdot)$ and the Volterra integral equation
    \begin{equation} \label{integral_equation}
        x(\tau)
        = a(\tau \mid t, w(\cdot))
        + \frac{1}{\mathrm{\Gamma}(\alpha)} \int_{t}^{\tau} \frac{f^\ast(\xi, x(\xi), \mu(\xi))}{(\tau - \xi)^{1 - \alpha}} \, \rd \xi
        \quad \forall \tau \in [t, T].
    \end{equation}

    For any $(t, w(\cdot)) \in G^0$ and any $u(\cdot) \in \mathcal{U}(t, T)$, taking into account that $f(\tau, x, u(\tau)) = f^\ast(\tau, x, \mu_{u(\cdot)}(\tau))$ for a.e. $\tau \in [t, T]$ and every $x \in \mathbb{R}^n$ according to \eqref{relaxed_control_via_usual_control} and \eqref{f_ast}, we have $x(\cdot \mid t, w(\cdot), u(\cdot)) = x(\cdot \mid t, w(\cdot), \mu_{u(\cdot)} (\cdot))$.
    Therefore, by definition \eqref{value_functional_G^0} of the value functional $\rho$ of the optimal control problem \eqref{system}, \eqref{cost_functional}, and recalling definition \eqref{M_0} of the set $\mathcal{M}_\ast(t, T)$, we get
    \begin{equation} \label{value_functional_M_0}
        \rho(t, w(\cdot))
        = \inf_{\mu(\cdot) \in \mathcal{M}_\ast(t, T)} \sigma\bigl( x(T \mid t, w(\cdot), \mu(\cdot)) \bigr)
        \quad \forall (t, w(\cdot)) \in G^0.
    \end{equation}

    Let us prove that the dependence of a motion $x(\cdot \mid t, w(\cdot), \mu(\cdot))$ of system \eqref{system} on the relaxed control $\mu(\cdot) \in \mathcal{M}(t, T)$ is continuous.
    \begin{proposition} \label{proposition_continuity_x_mu}
        Under assumptions $(f.1)$--$(f.3)$, the following mapping is continuous for every $(t, w(\cdot)) \in G^0$:
        \begin{equation*}
            \mathcal{M}(t, T) \ni \mu(\cdot) \mapsto x(\cdot \mid t, w(\cdot), \mu(\cdot)) \in \AC^\alpha([0, T], \mathbb{R}^n).
        \end{equation*}
    \end{proposition}
    {\it Proof}
        Let $\{\mu_i(\cdot)\}_{i \in \mathbb{N}_0} \subset \mathcal{M}(t, T)$ be such that $\mu_i(\cdot) \to \mu_0(\cdot)$ as $i \to \infty$.
        Here and below, we denote $\mathbb{N}_0 \coloneqq \mathbb{N} \cup \{0\}$.
        For every $i \in \mathbb{N}_0$, let us consider the motion $x_i(\cdot) \coloneqq x(\cdot \mid t, w(\cdot), \mu_i(\cdot))$ of system \eqref{system}.
        It is required to show that $x_i(\cdot) \to x_0(\cdot)$ as $i \to \infty$.

        By the sublinear growth property $(f^\ast.3)$, there exists $R \geq 0$ such that $\|x_i(\cdot)\|_{[0, T]} \leq R$ for all $i \in \mathbb{N}_0$ (see, e.g., \cite[Proposition 7.1, item (i)]{Gomoyunov_2020_SIAM}).
        Using the local Lipschitz continuity property $(f^\ast.2)$, let us take the corresponding number $\lambda_f \geq 0$.
        Based on \eqref{integral_equation}, for any $i \in \mathbb{N}$ and any $\tau \in [t, T]$, we derive
        \begin{align}
            & \|x_i(\tau) - x_0(\tau)\|
            \leq \frac{\lambda_f}{\mathrm{\Gamma}(\alpha)} \int_{t}^{\tau} \frac{\|x_i(\xi) - x_0(\xi)\|}{(\tau - \xi)^{1 - \alpha}} \, \rd \xi
            \nonumber \\
            & \ + \biggl\| \frac{1}{\mathrm{\Gamma}(\alpha)} \int_{t}^{\tau} \frac{f^\ast(\xi, x_0(\xi), \mu_i(\xi))}{(\tau - \xi)^{1 - \alpha}} \, \rd \xi
            - \frac{1}{\mathrm{\Gamma}(\alpha)} \int_{t}^{\tau} \frac{f^\ast(\xi, x_0(\xi), \mu_0(\xi))}{(\tau - \xi)^{1 - \alpha}} \, \rd \xi \biggr\|.
            \label{proposition_continuity_x_mu_proof_estimate}
        \end{align}

        For every $i \in \mathbb{N}_0$, let us define the function $h_i \colon [t, T] \to \mathbb{R}^n$ by (see \eqref{f_ast})
        \begin{align}
            h_i(\tau)
            & \coloneqq \frac{1}{\mathrm{\Gamma}(\alpha)} \int_{t}^{\tau} \frac{f^\ast(\xi, x_0(\xi), \mu_i(\xi))}{(\tau - \xi)^{1 - \alpha}} \, \rd \xi
            \nonumber \\
            & = \int_{t}^{\tau} \int_{P} \frac{f(\xi, x_0(\xi), u)}{\mathrm{\Gamma}(\alpha) (\tau - \xi)^{1 - \alpha}} \mu_i(\xi, \rd u) \, \rd \xi
            \quad \forall \tau \in [t, T].
            \label{proposition_continuity_x_mu_proof_notation}
        \end{align}
        Let us fix $\tau \in (t, T]$ and introduce the function $\mathfrak{b}^{[\tau]} \colon [t, T] \times P \to \mathbb{R}^n$ as follows:
        \begin{equation*}
            \mathfrak{b}^{[\tau]}(\xi, u)
            \coloneqq \begin{cases}
                \displaystyle
                \frac{f(\xi, x_0(\xi), u)}{\mathrm{\Gamma}(\alpha) (\tau - \xi)^{1 - \alpha}}, & \mbox{if } \xi \in [t, \tau), \\
                0, & \mbox{if } \xi \in [\tau, T],
              \end{cases}
        \end{equation*}
        for all $u \in P$.
        Note that the function $\mathfrak{b}^{[\tau]}$ (more precisely, every corresponding coordinate function) is a Carath\'{e}odory function in view of the continuity assumption $(f.1)$.
        Hence, recalling definition \eqref{convergence_M} of convergence in the space $\mathcal{M}(t, T)$, we get
        \begin{equation*}
            h_i(\tau)
            = \int_{t}^{T} \int_{P} \mathfrak{b}^{[\tau]}(\xi, u) \mu_i(\xi, \rd u) \, \rd \xi
            \to \int_{t}^{T} \int_{P} \mathfrak{b}^{[\tau]}(\xi, u) \mu_0(\tau, \rd u) \, \rd \xi
            = h_0(\tau)
        \end{equation*}
        as $i \to \infty$.
        Consequently, and since $h_i(t) = 0$ for all $i \in \mathbb{N}_0$, we conclude that $h_i(\tau) \to h_0(\tau)$ as $i \to \infty$ for all $\tau \in [t, T]$.
        In addition, for every $i \in \mathbb{N}_0$, we derive $\|f^\ast(\xi, x_0(\xi), \mu_i(\xi))\| \leq c_f(1 + R)$ for a.e. $\xi \in [t, T]$, where $c_f \geq 0$ is the number from $(f^\ast.3)$, and, therefore, the estimate below is valid (see, e.g., \cite[Proposition 2.1]{Gomoyunov_2019_FCAA_2}):
        \begin{equation*}
            \|h_i(\tau) - h_i(\tau^\prime)\|
            \leq \frac{2 c_f(1 + R)}{\mathrm{\Gamma}(\alpha + 1)} |\tau - \tau^\prime|^\alpha
            \quad \forall \tau, \tau^\prime \in [t, T].
        \end{equation*}
        Thus, the functions $h_i(\cdot)$ for all $i \in \mathbb{N}_0$ are equicontinuous, which together with the pointwise convergence established above implies that $\{h_i(\cdot)\}_{i \in \mathbb{N}}$ converges to $h_0(\cdot)$ as $i \to \infty$ uniformly on $[t, T]$ (see, e.g., \cite[Theorem I.5.3]{Warga_1972}).

        Now, let $\varepsilon > 0$ be given.
        Let us choose $\varepsilon_\ast > 0$ such that $\varepsilon_\ast \mathrm{E}_\alpha(\lambda_f T^\alpha) \leq \varepsilon$, where $\mathrm{E}_\alpha$ is the Mittag-Leffler function (see, e.g., \cite[Chapter 3]{Gorenflo_Kilbas_Mainardi_Rogosin_2014}).
        Then, there exists $i_\ast \in \mathbb{N}$ such that $\|h_i(\tau) - h_0(\tau)\| \leq \varepsilon_\ast$ for all $i \in \mathbb{N}$ with $i \geq i_\ast$ and all $\tau \in [t, T]$.
        Let $i \in \mathbb{N}$ with $i \geq i_\ast$ be fixed.
        By virtue of \eqref{proposition_continuity_x_mu_proof_estimate} and \eqref{proposition_continuity_x_mu_proof_notation}, we have
        \begin{equation*}
            \|x_i(\tau) - x_0(\tau)\|
            \leq \frac{\lambda_f}{\mathrm{\Gamma}(\alpha)} \int_{t}^{\tau} \frac{\|x_i(\xi) - x_0(\xi)\|}{(\tau - \xi)^{1 - \alpha}} \, \rd \xi
            + \varepsilon_\ast
            \quad \forall \tau \in [t, T],
        \end{equation*}
        wherefrom, applying a Gronwall type inequality (see, e.g., \cite[Lemma 6.19]{Diethelm_2010}), we derive
        \begin{equation*}
            \|x_i(\tau) - x_0(\tau)\|
            \leq \varepsilon_\ast \mathrm{E}_\alpha(\lambda_f (\tau - t)^\alpha)
            \leq \varepsilon
            \quad \forall \tau \in [t, T].
        \end{equation*}
        Hence, taking into account that $x_i(\tau) = w(\tau) = x_0(\tau)$ for all $\tau \in [0, t]$, we arrive at the estimate $\|x_i(\cdot) - x_0(\cdot)\|_{[0, T]} \leq \varepsilon$.
        Thus, $x_i(\cdot) \to x_0(\cdot)$ as $i \to \infty$, and the proof is complete.
    \qed

    In particular, due to Proposition \ref{proposition_continuity_x_mu}, continuity of the function $\sigma$ (see assumption $(\sigma.1)$), compactness of the space $\mathcal{M} (t, T)$, and the fact that $\mathcal{M}_\ast(t, T)$ is dense in $\mathcal{M}(t, T)$, we have
    \begin{equation*}
        \inf_{\mu(\cdot) \in \mathcal{M}_\ast(t, T)} \sigma \bigl( x(T \mid t, w(\cdot), \mu(\cdot)) \bigr)
        = \min_{\mu(\cdot) \in \mathcal{M} (t, T)} \sigma\bigl( x(T \mid t, w(\cdot), \mu(\cdot)) \bigr)
    \end{equation*}
    for all $(t, w(\cdot)) \in G^0$.
    As a result (see \eqref{value_functional_M_0}), we get the following representation of the value functional $\rho$ of the optimal control problem \eqref{system}, \eqref{cost_functional}:
    \begin{equation} \label{value_functional_mu}
        \rho(t, w(\cdot))
        = \min_{\mu(\cdot) \in \mathcal{M} (t, T)} \sigma\bigl( x(T \mid t, w(\cdot), \mu(\cdot)) \bigr)
        \quad \forall (t, w(\cdot)) \in G^0.
    \end{equation}

    Another difficulty that arises in applying Theorem \ref{theorem_envelope_theorem} to our setting is related to the fact that the space $\mathcal{M}(t, T)$, over which the minimum is taken in \eqref{value_functional_mu}, depends on the time variable $t$, while the space of parameters $\mathcal{L}$ in \eqref{varphi_is_envelope} does not.
    We overcome this difficulty by making for every fixed initial position $(t, w(\cdot)) \in G^0$ a change of the time variable in the dynamic equation \eqref{system} that allows us to obtain the unified control interval $[0, 1]$ instead of $[t, T]$.

\section{Change of Time Variable}

    Let $(t, w(\cdot)) \in G^0$ be fixed.
    Let us denote $\mathcal{N} \coloneqq \mathcal{M}(0, 1)$, where the metric space $\mathcal{M}(0, 1)$ is defined in a similar way as $\mathcal{M}(t, T)$ in Section \ref{section_relaxed_controls}.
    Let us introduce the mapping $\pi \colon \mathcal{M}(t, T) \to \mathcal{N}$ that assigns to each function $\mu(\cdot) \in \mathcal{M}(t, T)$ the function $\nu(\cdot) \in \mathcal{N}$ such that
    \begin{equation} \label{pi}
        \nu(\vartheta)
        \coloneqq \mu(t + \vartheta (T - t))
        \quad \text{for a.e. } \vartheta \in [0, 1].
    \end{equation}
    Note that the mapping $\pi$ is well-defined and one-to-one and its inverse $\pi^{- 1}$ assigns to each function $\nu(\cdot) \in \mathcal{N}$ the function $\mu(\cdot) \in \mathcal{M}(t, T)$ such that $\mu(\tau) \coloneqq \nu((\tau - t) / (T - t))$ for a.e. $\tau \in [t, T]$.
    In addition, it follows directly from the definition of convergence in the spaces $\mathcal{M}(t, T)$ and $\mathcal{N}$ (see \eqref{convergence_M}) that the mappings $\pi$ and $\pi^{- 1}$ are continuous.

    For every $\nu(\cdot) \in \mathcal{N}$, let us consider the {\it auxiliary Volterra integral equation}
    \begin{align}
        y(\vartheta)
        & = a(t + \vartheta (T - t) \mid t, w(\cdot))
        \nonumber \\
        & \quad + \frac{(T - t)^\alpha}{\mathrm{\Gamma}(\alpha)} \int_{0}^{\vartheta}
        \frac{f^\ast(t + \zeta (T - t), y(\zeta), \nu(\zeta))}{(\vartheta - \zeta)^{1 - \alpha}} \, \rd \zeta
        \quad \forall \vartheta \in [0, 1],
        \label{y_integral_equation}
    \end{align}
    where the function $a(\cdot \mid t, w(\cdot))$ is defined according to \eqref{a} and the function $f^\ast$ is given by \eqref{f_ast}.
    By a {\it solution} of the integral equation \eqref{y_integral_equation}, we mean a function $y(\cdot) \in \Cont([0, 1], \mathbb{R}^n)$ that satisfies this equation.

    Similarly to the proof of \cite[Lemma 4.1]{Gomoyunov_2022_FCAA}, making the change of variables $\vartheta \coloneqq (\tau - t) / (T - t)$ in the (original) integral equation \eqref{integral_equation}, we obtain
    \begin{proposition} \label{proposition_y_x}
        Under assumptions $(f.1)$--$(f.3)$, for any $(t, w(\cdot)) \in G^0$ and any $\nu(\cdot) \in \mathcal{N}$, there exists a unique solution $y(\cdot) \coloneqq y(\cdot \mid t, w(\cdot), \nu(\cdot))$ of the integral equation \eqref{y_integral_equation}.
        Moreover, the equality
        \begin{equation*}
            y(\vartheta)
            = x(t + \vartheta (T - t))
            \quad \forall \vartheta \in [0, 1]
        \end{equation*}
        holds, where $x(\cdot) \coloneqq x(\cdot \mid t, w(\cdot), \mu(\cdot))$ is the motion of system \eqref{system} generated from the initial position $(t, w(\cdot))$ by the relaxed control $\mu(\cdot) \coloneqq \pi^{- 1}(\nu(\cdot))$.
    \end{proposition}

    Based on relation \eqref{value_functional_mu} and Proposition \ref{proposition_y_x}, we derive the following representation of the value functional $\rho$ of the optimal control problem \eqref{system}, \eqref{cost_functional}:
    \begin{equation} \label{value_functional_representation}
        \rho(t, w(\cdot))
        = \min_{\nu(\cdot) \in \mathcal{N}} \sigma\bigl( y(1 \mid t, w(\cdot), \nu(\cdot)) \bigr)
        \quad \forall (t, w(\cdot)) \in G^0.
    \end{equation}
    We see that this representation is of the form \eqref{varphi_is_envelope} with $\mathcal{L} \coloneqq \mathcal{N}$ and
    \begin{equation} \label{psi}
        \psi(t, w(\cdot), \nu(\cdot))
        \coloneqq \sigma \bigl( y(1 \mid t, w(\cdot), \nu(\cdot)) \bigr)
        \quad \forall (t, w(\cdot)) \in G^0 \quad \forall \nu(\cdot) \in \mathcal{N}.
    \end{equation}
    Thus, in order to apply Theorem \ref{theorem_envelope_theorem}, we need to establish the required continuity and directional differentiability properties of this functional $\psi$, which mainly reduces to studying the corresponding properties of the mapping
    \begin{equation*}
        G^0 \times \mathcal{N} \ni ((t, w(\cdot)), \nu(\cdot)) \mapsto y(1 \mid t, w(\cdot), \nu(\cdot)) \in \mathbb{R}^n.
    \end{equation*}

\section{Continuity Properties}
\label{section_continuity}

    The goal of this section is to prove
    \begin{lemma} \label{lemma_y_continuous}
        Let assumptions $(f.1)$--$(f.3)$ hold.
        Then, the following mapping is continuous:
        \begin{equation} \label{lemma_y_continuous_main}
            G^0 \times \mathcal{N} \ni ((t, w(\cdot)), \nu(\cdot))
            \mapsto y(\cdot \mid t, w(\cdot), \nu(\cdot)) \in \Cont([0, 1], \mathbb{R}^n),
        \end{equation}
        where $y(\cdot \mid t, w(\cdot), \nu(\cdot))$ is the solution of the integral equation \eqref{y_integral_equation}.
    \end{lemma}
    {\it Proof}
        1.
            We first note that, for every fixed $(t, w(\cdot)) \in G^0$, the mapping
            \begin{equation*}
                \mathcal{N} \ni \nu(\cdot) \mapsto y(\cdot \mid t, w(\cdot), \nu(\cdot)) \in \Cont([0, 1], \mathbb{R}^n)
            \end{equation*}
            is continuous.
            Indeed, let $\{\nu_i(\cdot)\}_{i \in \mathbb{N}_0} \subset \mathcal{N}$ be such that $\nu_i(\cdot) \to \nu_0(\cdot)$ as $i \to \infty$.
            For every $i \in \mathbb{N}_0$, let us consider the solution $y_i(\cdot) \coloneqq y(\cdot \mid t, w(\cdot), \nu_i(\cdot))$ of the integral equation \eqref{y_integral_equation} and the motion $x_i(\cdot) \coloneqq x(\cdot \mid t, w(\cdot), \mu_i(\cdot))$ of system \eqref{system} generated by the relaxed control $\mu_i(\cdot) \coloneqq \pi^{- 1}(\nu_i(\cdot))$ (see \eqref{pi}).
            Since the mapping $\pi^{- 1}$ is continuous, we have $\mu_i(\cdot) \to \mu_0(\cdot)$ as $i \to \infty$, and, therefore, $x_i(\cdot) \to x_0(\cdot)$ as $i \to \infty$ by Proposition \ref{proposition_continuity_x_mu}.
            Consequently, taking into account that $y_i(\vartheta) = x_i(t + \vartheta (T - t))$ for all $\vartheta \in [0, 1]$ and all $i \in \mathbb{N}_0$ due to Proposition \ref{proposition_y_x}, we conclude that $y_i(\cdot) \to y_0(\cdot)$ as $i \to \infty$.

        2.
            Let us take an arbitrary compact set $K \subset G^0$ and prove that, for every $\varepsilon > 0$, there exists $\delta > 0$ such that, for any $(t, w(\cdot))$, $(t^\prime, w^\prime(\cdot)) \in K$ satisfying the condition $\dist((t, w(\cdot)), (t^\prime, w^\prime(\cdot))) \leq \delta$ (see \eqref{dist}) and any $\nu(\cdot) \in \mathcal{N}$,
            \begin{equation*}
                \|y(\cdot \mid t, w(\cdot), \nu(\cdot)) - y(\cdot \mid t^\prime, w^\prime(\cdot), \nu(\cdot))\|_{[0, 1]}
                \leq \varepsilon.
            \end{equation*}

            In view of property $(f^\ast.3)$, there exists $R \geq 0$ such that the inequality $\|x(\cdot \mid t, w(\cdot), \mu(\cdot))\|_{[0, T]} \leq R$ is valid for all $(t, w(\cdot)) \in K$ and all $\mu(\cdot) \in \mathcal{M}(t, T)$ (see, e.g., \cite[Proposition 7.1, item (i)]{Gomoyunov_2020_SIAM}).
            In particular, by Proposition \ref{proposition_y_x}, we obtain $\|y(\cdot \mid t, w(\cdot), \nu(\cdot))\|_{[0, 1]} \leq R$ for all $(t, w(\cdot)) \in K$ and all $\nu(\cdot) \in \mathcal{N}$.
            Using property $(f^\ast.2)$, let us take the number $\lambda_f \geq 0$ that corresponds to $R$.

            Let $\varepsilon > 0$ be given.
            Let us choose $\varepsilon_\ast > 0$ such that $4 \varepsilon_\ast \mathrm{E}_\alpha (T^\alpha \lambda_f) \leq \varepsilon$, where $\mathrm{E}_\alpha$ denotes again the Mittag-Leffler function.
            Due to property $(f^\ast.1)$, there exists $\delta_1 > 0$ such that, for any $\tau$, $\tau^\prime \in [0, T]$ with $|\tau - \tau^\prime| \leq \delta_1$, any $x \in B(R)$, and any $\mu \in \rpm(P)$,
            \begin{equation*}
                \frac{T^\alpha \|f^\ast(\tau, x, \mu) - f^\ast(\tau^\prime, x, \mu)\|}{\mathrm{\Gamma}(\alpha + 1)}
                \leq \varepsilon_\ast.
            \end{equation*}
            Further, according to \cite[Lemma 3]{Gomoyunov_2020_DE}, the mapping
            \begin{equation} \label{a_mapping}
                G^0 \ni (t, w(\cdot)) \mapsto a(\cdot \mid t, w(\cdot)) \in \AC^\alpha([0, T], \mathbb{R}^n)
            \end{equation}
            is continuous, where the function $a(\cdot \mid t, w(\cdot))$ is defined by \eqref{a}.
            Hence, owing to compactness of the set $K$, we conclude that there exists $\delta_2 > 0$ such that
            \begin{equation*}
                \| a(\cdot \mid t, w(\cdot)) - a(\cdot \mid t^\prime, w^\prime(\cdot)) \|_{[0, T]}
                \leq \varepsilon_\ast
            \end{equation*}
            for all $(t, w(\cdot))$, $(t^\prime, w^\prime(\cdot)) \in K$ with $\dist((t, w(\cdot)), (t^\prime, w^\prime(\cdot))) \leq \delta_2$.
            In addition, we obtain that the set $\{a(\cdot \mid t, w(\cdot)) \in \AC^\alpha([0, T], \mathbb{R}^n) \colon (t, w(\cdot)) \in K\}$ is compact, which implies that the functions $a(\cdot \mid t, w(\cdot))$ with $(t, w(\cdot)) \in K$ are equicontinuous.
            Then, there exists $\delta_3 > 0$ such that, for any $(t, w(\cdot)) \in K$ and any $\tau$, $\tau^\prime \in [0, T]$ with $|\tau - \tau^\prime| \leq \delta_3$,
            \begin{equation*}
                \|a(\tau \mid t, w(\cdot)) - a(\tau^\prime \mid t, w(\cdot))\|
                \leq \varepsilon_\ast.
            \end{equation*}
            Finally, let us take $\delta_4 > 0$ satisfying the condition
            \begin{equation*}
                \frac{\delta_4^\alpha c_f (1 + R)}{\mathrm{\Gamma}(\alpha + 1)}
                \leq \varepsilon_\ast,
            \end{equation*}
            where $c_f \geq 0$ is the number from property $(f^\ast.3)$.
            Let us verify that the claim holds for $\delta \coloneqq \min_{i \in \overline{1, 4}} \delta_i > 0$.

            Let $(t, w(\cdot))$, $(t^\prime, w^\prime(\cdot)) \in K$ with $\dist((t, w(\cdot)), (t^\prime, w^\prime(\cdot))) \leq \delta$ and $\nu(\cdot) \in \mathcal{N}$ be fixed.
            Let us denote $y(\cdot) \coloneqq y(\cdot \mid t, w(\cdot), \nu(\cdot))$ and $y^\prime(\cdot) \coloneqq y(\cdot \mid t^\prime, w^\prime(\cdot), \nu(\cdot))$.
            According to \eqref{y_integral_equation}, for any $\vartheta \in [0, 1]$, we have
            \begin{align*}
                & \|y(\vartheta) - y^\prime(\vartheta)\|
                \leq \|a(t + \vartheta (T - t) \mid t, w(\cdot)) - a(t^\prime + \vartheta (T - t^\prime) \mid t^\prime, w^\prime(\cdot))\| \\
                & \quad + \frac{|(T - t)^\alpha - (T - t^\prime)^\alpha|}{\mathrm{\Gamma}(\alpha)}
                \int_{0}^{\vartheta} \frac{\|f^\ast(t + \zeta (T - t), y(\zeta), \nu(\zeta))\|}{(\vartheta - \zeta)^{1 - \alpha}} \, \rd \zeta \\
                & \quad + \frac{(T - t^\prime)^\alpha}{\mathrm{\Gamma}(\alpha)}
                \int_{0}^{\vartheta} \biggl\| \frac{f^\ast(t + \zeta (T - t), y(\zeta), \nu(\zeta))}{(\vartheta - \zeta)^{1 - \alpha}} \\
                & \quad - \frac{f^\ast(t^\prime + \zeta (T - t^\prime), y^\prime(\zeta), \nu(\zeta))}{(\vartheta - \zeta)^{1 - \alpha}} \biggr\| \, \rd \zeta.
            \end{align*}
            Let us estimate each of the terms from the right-hand side of this inequality separately.
            For the first term, denoted by $r^{(1)}$, we derive
            \begin{align*}
                r^{(1)}
                & \leq \| a(t + \vartheta (T - t) \mid t, w(\cdot)) - a(t^\prime + \vartheta (T - t^\prime) \mid t, w(\cdot)) \| \\
                & \quad + \| a(t^\prime + \vartheta (T - t^\prime) \mid t, w(\cdot)) - a(t^\prime + \vartheta (T - t^\prime) \mid t^\prime, w^\prime(\cdot)) \| \\
                & \leq 2 \varepsilon_\ast.
            \end{align*}
            For the second term, denoted by $r^{(2)}$, we obtain
            \begin{equation*}
                r^{(2)}
                \leq \frac{|t - t^\prime|^\alpha c_f (1 + R)}{\mathrm{\Gamma}(\alpha)}
                \int_{0}^{\vartheta} \frac{\rd \zeta}{(\vartheta - \zeta)^{1 - \alpha}}
                \leq \frac{\delta^\alpha c_f (1 + R)}{\mathrm{\Gamma}(\alpha + 1)}
                \leq \varepsilon_\ast.
            \end{equation*}
            For the third term, denoted by $r^{(3)}$, we get
            \begin{align*}
                & r^{(3)}
                \leq \frac{T^\alpha \lambda_f}{\mathrm{\Gamma}(\alpha)} \int_{0}^{\vartheta}
                \frac{\|y(\zeta) - y^\prime(\zeta)\|}{(\vartheta - \zeta)^{1 - \alpha}} \, \rd \zeta \\
                & \ \ + \frac{T^\alpha}{\mathrm{\Gamma}(\alpha)} \int_{0}^{\vartheta}
                \frac{\|f^\ast(t + \zeta (T - t), y(\zeta), \nu(\zeta)) - f^\ast(t^\prime + \zeta (T - t^\prime), y(\zeta), \nu(\zeta))\|}
                {(\vartheta - \zeta)^{1 - \alpha}} \, \rd \zeta \\
                & \ \ \leq \frac{T^\alpha \lambda_f}{\mathrm{\Gamma}(\alpha)} \int_{0}^{\vartheta}
                \frac{\|y(\zeta) - y^\prime(\zeta)\|}{(\vartheta - \zeta)^{1 - \alpha}} \, \rd \zeta
                + \varepsilon_\ast.
            \end{align*}
            As a result, we arrive at the estimate
            \begin{equation*}
                \| y(\vartheta) - y^\prime(\vartheta) \|
                \leq 4 \varepsilon_\ast
                + \frac{T^\alpha \lambda_f}{\mathrm{\Gamma}(\alpha)} \int_{0}^{\vartheta} \frac{\|y(\zeta) - y^\prime(\zeta)\|}{(\vartheta - \zeta)^{1 - \alpha}} \, \rd \zeta
                \quad \forall \vartheta \in [0, 1].
            \end{equation*}
            Then, applying a Gronwall type inequality (see, e.g., \cite[Lemma 6.19]{Diethelm_2010}), we conclude that
            \begin{equation*}
                \| y(\vartheta) - y^\prime(\vartheta) \|
                \leq 4 \varepsilon_\ast \mathrm{E}_\alpha (T^\alpha \lambda_f \vartheta^\alpha)
                \leq \varepsilon
                \quad \forall \vartheta \in [0, 1],
            \end{equation*}
            which completes the proof of the claim.

        3.
            At the final step, let us show that mapping \eqref{lemma_y_continuous_main} is continuous.
            Let $\{((t_i, w_i(\cdot)), \nu_i(\cdot))\}_{i \in \mathbb{N}_0} \subset G^0 \times \mathcal{N}$ and $((t_i, w_i(\cdot)), \nu_i(\cdot)) \to ((t_0, w_0(\cdot)), \nu_0(\cdot))$ as $i \to \infty$.
            Then, we have $\|y(\cdot \mid t_0, w_0(\cdot), \nu_i(\cdot)) - y(\cdot \mid t_0, w_0(\cdot), \nu_0(\cdot))\|_{[0, 1]} \to 0$ as $i \to \infty$ by the first part of the proof.
            Due to the second part of the proof, considering the compact set $K \coloneqq \{(t_i, w_i(\cdot)) \in G^0 \colon i \in \mathbb{N}_0\}$, we obtain that $\|y(\cdot \mid t_i, w_i(\cdot), \nu_i(\cdot)) - y(\cdot \mid t_0, w_0(\cdot), \nu_i(\cdot))\|_{[0, 1]} \to 0$ as $i \to \infty$.
            Hence, we get the desired convergence $\|y(\cdot \mid t_i, w_i(\cdot), \nu_i(\cdot)) - y(\cdot \mid t_0, w_0(\cdot), \nu_0(\cdot))\|_{[0, 1]} \to 0$ as $i \to \infty$.
            The lemma is proved.
    \qed

    Lemma \ref{lemma_y_continuous} immediately implies
    \begin{corollary} \label{corollary_psi_is_continuous}
        Let assumptions $(f.1)$--$(f.3)$ and $(\sigma.1)$ hold.
        Then, the functional $\psi$ from \eqref{psi} is continuous.
    \end{corollary}

\section{Directional Differentiability Properties}
\label{section_differentiability}

    \subsection{Auxiliary Notation}
    \label{subsection_auxiliary_notation}

        Let $N \geq c_f$ be given, where $c_f \geq 0$ is the number from property $(f^\ast.3)$.
        Let us denote by $X_N$ the set of all functions $x(\cdot) \in \AC^\alpha([0, T], \mathbb{R}^n)$ such that $\|x(0)\| \leq N$ and $\|(^C D^\alpha x)(\tau)\| \leq N (1 + \|x(\tau)\|)$ for a.e. $\tau \in [0, T]$.
        The set $X_N$ is non-empty and compact (see, e.g., \cite[Theorems 1 and 2]{Gomoyunov_2020_DE}).
        In particular, there exists $R_N \geq 0$ such that, for any $x(\cdot) \in X_N$,
        \begin{equation} \label{R_N}
            \|x(\cdot)\|_{[0, T]}
            \leq R_N,
            \quad \|(^C D^\alpha x)(\tau)\|
            \leq R_N
            \quad \text{for a.e. } \tau \in [0, T].
        \end{equation}
        Further, let us consider the set $G_N$ of all points $(t, w(\cdot)) \in G$ for each of which there exists $x(\cdot) \in X_N$ such that $w(\cdot) = x_t(\cdot)$ (see \eqref{x_t}).
        The set $G_N$ is non-empty and compact (see, e.g., the proof of \cite[Proposition 10.2, item (iv)]{Gomoyunov_2020_SIAM}).
        Let us note that, due to property $(f^\ast.3)$, for any $(t, w(\cdot)) \in G_N \cap G^0$ and any $\mu(\cdot) \in \mathcal{M}(t, T)$, the motion $x(\cdot) \coloneqq x(\cdot \mid t, w(\cdot), \mu(\cdot))$ of system \eqref{system} satisfies the inclusion $x(\cdot) \in X_N$ (see, e.g., \cite[Theorem 3]{Gomoyunov_2020_DE}), and, consequently, $(\tau, x_\tau(\cdot)) \in G_N$ for all $\tau \in [0, T]$.
        In addition, for any $(t, w(\cdot)) \in G_N \cap G^0$, any $M \geq 0$, and any $f \in B(M)$, the function $x^{(f)}(\cdot) \coloneqq x^{(f)}(\cdot \mid t, w(\cdot))$ defined by \eqref{x^f} satisfies the inclusions $x^{(f)}(\cdot) \in X_{N + M}$ and $(\tau, x^{(f)}_\tau(\cdot)) \in G_{N + M}$ for all $\tau \in [0, T]$.
        Finally, for every $\eta \in (0, T)$, let us introduce the non-empty and compact set
        \begin{equation*}
            G_N^\eta
            \coloneqq \bigl\{ (t, w(\cdot)) \in G_N \colon
            t \leq T - \eta \bigr\}.
        \end{equation*}
        It is clear that the union of $X_N$ for all $N \geq c_f$ gives $\AC^\alpha([0, T], \mathbb{R}^n)$, and, therefore, the union of $G_N^\eta$ for all $N \geq c_f$ and all $\eta \in (0, T)$ coincides with $G^0$.

        In particular, from Proposition \ref{proposition_y_x} and Lemma \ref{lemma_y_continuous}, we derive
        \begin{corollary} \label{corollary_y_uniform_continuity}
            Let assumptions $(f.1)$--$(f.3)$ hold, and let $N \geq c_f$, $M \geq 0$, and $\eta \in (0, T)$.
            Then, for any $(t, w(\cdot)) \in G_N^\eta$, any $f \in B(M)$, any $\nu(\cdot) \in \mathcal{N}$, and any $\delta \in [0, \eta / 2]$, the inequality below is valid:
            \begin{equation} \label{corollary_y_uniform_continuity_R_N+M}
                \|y(\cdot \mid t + \delta, x^{(f)}_{t + \delta} (\cdot \mid t, w(\cdot)), \nu(\cdot))\|_{[0, 1]}
                \leq R_{N + M},
            \end{equation}
            where $y(\cdot \mid t + \delta, x^{(f)}_{t + \delta} (\cdot \mid t, w(\cdot)), \nu(\cdot))$ is the solution of the integral equation \eqref{y_integral_equation}.
            In addition,
            \begin{equation} \label{corollary_y_uniform_continuity_main}
                \| y(\cdot \mid t + \delta, x^{(f)}_{t + \delta}(\cdot \mid t, w(\cdot)), \nu(\cdot))
                - y(\cdot \mid t, w(\cdot), \nu(\cdot))\|_{[0, 1]}
                \to 0
            \end{equation}
            as $\delta \to 0^+$ uniformly in $(t, w(\cdot)) \in G_N^\eta$, $f \in B(M)$, and $\nu(\cdot) \in \mathcal{N}$.
        \end{corollary}
        {\it Proof}
        1.
            Let $(t, w(\cdot)) \in G_N^\eta$, $f \in B(M)$, $\nu(\cdot) \in \mathcal{N}$, and $\delta \in [0, \eta / 2]$ be fixed.
            Let us put $\mu(\cdot) \coloneqq \pi^{- 1}(\nu(\cdot))$ (see \eqref{pi}) and consider the corresponding motion $x(\cdot) \coloneqq x(\cdot \mid t + \delta, x^{(f)}_{t + \delta}(\cdot \mid t, w(\cdot)), \mu(\cdot))$ of system \eqref{system}.
            Taking into account that $(t + \delta, x^{(f)}_{t + \delta}(\cdot \mid t, w(\cdot))) \in G_{N + M} \cap G^0$, we obtain $x(\cdot) \in X_{N + M}$, and, hence, $\|x(\cdot)\|_{[0, T]} \leq R_{N + M}$.
            Then, applying Proposition \ref{proposition_y_x}, for any $\vartheta \in [0, 1]$, we get
            \begin{equation*}
                \| y(\vartheta \mid t + \delta, x^{(f)}_{t + \delta} (\cdot \mid t, w(\cdot)), \nu(\cdot)) \|
                = \|x(t + \delta + \vartheta (T - t - \delta))\|
                \leq R_{N + M}.
            \end{equation*}

        2.
            Due to continuity of mapping \eqref{a_mapping} and definition \eqref{x^f} of the function $x^{(f)}(\cdot \mid t, w(\cdot))$, the mapping
            \begin{equation*}
                G_N^\eta \times B(M) \ni ((t, w(\cdot)), f) \mapsto x^{(f)}(\cdot \mid t, w(\cdot)) \in \AC^\alpha([0, T], \mathbb{R}^n)
            \end{equation*}
            is continuous.
            Therefore, by continuity of mapping \eqref{general_mapping}, the mapping
            \begin{equation*}
                [0, \eta / 2] \times G_N^\eta \times B(M) \ni (\delta, (t, w(\cdot)), f)
                \mapsto (t + \delta, x^{(f)}_{t + \delta}(\cdot \mid t, w(\cdot))) \in G
            \end{equation*}
            is continuous.
            Consequently, it follows from Lemma \ref{lemma_y_continuous} that the mapping
            \begin{align*}
                & [0, \eta / 2] \times G_N^\eta \times B(M) \times \mathcal{N} \ni (\delta, (t, w(\cdot)), f, \nu(\cdot)) \\
                & \quad \mapsto y(\cdot \mid t + \delta, x^{(f)}_{t + \delta}(\cdot \mid t, w(\cdot)), \nu(\cdot))
                \in \Cont([0, 1], \mathbb{R}^n)
            \end{align*}
            is continuous.
            Then, for any $(t, w(\cdot)) \in G_N^\eta$, any $f \in B(M)$, and any $\nu(\cdot) \in \mathcal{N}$, since $x^{(f)}_t(\cdot \mid t, w(\cdot)) = w(\cdot)$ by \eqref{x^f}, we conclude that convergence \eqref{corollary_y_uniform_continuity_main} takes place.
            Moreover, recalling that the sets $G_N^\eta$, $B(M)$, and $\mathcal{N}$ are compact, we obtain that this convergence is uniform (see, e.g., \cite[Theorem I.2.15]{Warga_1972}).
        \qed

    \subsection{Directional Differentiability of Solutions of Auxiliary Integral Equation}

        From now on, let us assume that the function $f$ from the right-hand side of the dynamic equation \eqref{system} satisfies the following additional condition:
        \begin{description}
            \item[$(f.4)$]
                The partial derivatives
                \begin{equation*}
                    \frac{\partial f}{\partial \tau} \colon [0, T] \times \mathbb{R}^n \times P \to \mathbb{R}^n,
                    \quad \frac{\partial f}{\partial x} \colon [0, T] \times \mathbb{R}^n \times P \to \mathbb{R}^{n \times n}
                \end{equation*}
                exist and are continuous.
                In the cases $\tau = 0$ and $\tau = T$, the partial derivative $\partial f / \partial \tau$ is understood as the corresponding one-sided derivative.
        \end{description}
        Then, in particular, the function $f^\ast$ given by \eqref{f_ast} possesses the property below:
        \begin{description}
            \item[$(f^\ast.4)$]
                For any $\tau \in [0, T]$, any $x \in \mathbb{R}^n$, and any $\mu \in \rpm(P)$, the equalities
                \begin{equation} \label{f_ast_partial_derivatives}
                    \frac{\partial f^\ast}{\partial \tau}(\tau, x, \mu)
                    = \int_{P} \frac{\partial f}{\partial \tau}(\tau, x, u) \mu(\rd u),
                    \quad \frac{\partial f^\ast}{\partial x}(\tau, x, \mu)
                    = \int_{P} \frac{\partial f}{\partial x}(\tau, x, u) \mu(\rd u)
                \end{equation}
                hold, and, moreover, the mappings $\partial f^\ast / \partial \tau \colon [0, T] \times \mathbb{R}^n \times \rpm(P) \to \mathbb{R}^n$ and $\partial f^\ast / \partial x \colon [0, T] \times \mathbb{R}^n \times \rpm(P) \to \mathbb{R}^{n \times n}$ are continuous.
        \end{description}
        Let us note also that condition $(f.4)$ implies condition $(f.2)$.

        Let $\Cont^{1 - \alpha}((0, 1], \mathbb{R}^n)$ be the set of all continuous functions $q \colon (0, 1] \to \mathbb{R}^n$ such that the function $(0, 1] \ni \vartheta \mapsto \vartheta^{1 - \alpha} \|q(\vartheta)\|$ is bounded.
        In a similar way, let us introduce the set $\Cont^{1 - \alpha}((0, 1], \mathbb{R}^{n \times n})$.

        Let $(t, w(\cdot)) \in G^0$ be fixed.
        Let us denote
        \begin{align}
            q(\vartheta)
            & \coloneqq - \frac{(1 - \alpha) (1 - \vartheta)}{\mathrm{\Gamma}(\alpha)}
            \int_{0}^{t} \frac{(^C D^\alpha w)(\xi)}{(t + \vartheta (T - t) - \xi)^{2 - \alpha}} \, \rd \xi,
            \nonumber \\
            Q(\vartheta)
            & \coloneqq \frac{\Id}{\mathrm{\Gamma}(\alpha) \vartheta^{1 - \alpha} (T - t)^{1 - \alpha}}
            \quad \forall \vartheta \in (0, 1],
            \label{r_R}
        \end{align}
        where $\Id \in \mathbb{R}^{n \times n}$ is the identity matrix.
        Note that $q(\cdot) \in \Cont^{1 - \alpha}((0, 1], \mathbb{R}^n)$ and $Q(\cdot) \in \Cont^{1 - \alpha}((0, 1], \mathbb{R}^{n \times n})$ (see, e.g., \cite[Lemma 3]{Gomoyunov_2022_FCAA}).
        Further, let $\nu(\cdot) \in \mathcal{N}$ be also fixed.
        Let us take the corresponding solution $y(\cdot) \coloneqq y(\cdot \mid t, w(\cdot), \nu(\cdot))$ of integral equation \eqref{y_integral_equation}, denote
        \begin{align}
            A(\vartheta, u)
            & \coloneqq (T - t)^\alpha \frac{\partial f}{\partial x} (t + \vartheta (T - t), y(\vartheta), u),
            \nonumber \\
            b(\vartheta, u)
            & \coloneqq (1 - \vartheta) (T - t)^\alpha \frac{\partial f}{\partial \tau} (t + \vartheta (T - t), y(\vartheta), u)
            \nonumber \\
            & \quad - \frac{\alpha f(t + \vartheta (T - t), y(\vartheta), u)}{(T - t)^{1 - \alpha}}
            \quad \forall \vartheta \in [0, 1] \quad \forall u \in P,
            \label{A_b}
        \end{align}
        and, accordingly, for any $\vartheta \in [0, 1]$ and any $\mu \in \rpm(P)$, put
        \begin{equation} \label{A_b^ast}
            A^\ast(\vartheta, \mu)
            \coloneqq \int_{P} A(\vartheta, u) \mu(\rd u),
            \quad b^\ast(\vartheta, \mu)
            \coloneqq \int_{P} b(\vartheta, u) \mu(\rd u).
        \end{equation}
        The functions $A \colon [0, 1] \times P \to \mathbb{R}^{n \times n}$ and $b \colon [0, 1] \times P \to \mathbb{R}^n$ are continuous in view of assumptions $(f.1)$ and $(f.4)$, which, in particular, implies that the functions $A^\ast \colon [0, 1] \times \rpm(P) \to \mathbb{R}^{n \times n}$ and $b^\ast \colon [0, 1] \times \rpm(P) \to \mathbb{R}^n$ are well-defined and continuous.
        Let us consider two {\it linear Volterra integral equations}
        \begin{equation} \label{z^1_integral_equation}
            z(\vartheta)
            = q(\vartheta) + \frac{1}{\mathrm{\Gamma}(\alpha)} \int_{0}^{\vartheta}
            \frac{A^\ast(\zeta, \nu(\zeta)) z(\zeta) + b^\ast(\zeta, \nu(\zeta))}{(\vartheta - \zeta)^{1 - \alpha}} \, \rd \zeta
            \quad \forall \vartheta \in (0, 1]
        \end{equation}
        and
        \begin{equation} \label{z^2_integral_equation}
            Z(\vartheta)
            = Q(\vartheta) + \frac{1}{\mathrm{\Gamma}(\alpha)} \int_{0}^{\vartheta}
            \frac{A^\ast(\zeta, \nu(\zeta)) Z(\zeta)}{(\vartheta - \zeta)^{1 - \alpha}} \, \rd \zeta
            \quad \forall \vartheta \in (0, 1].
        \end{equation}
        By a {\it solution} of equation \eqref{z^1_integral_equation}, we mean a function $z(\cdot) \in \Cont^{1 - \alpha}((0, 1], \mathbb{R}^n)$ satisfying this equation.
        Arguing similarly to the proof of \cite[Proposition 1]{Gomoyunov_2022_FCAA}, it can be shown that such a solution $z(\cdot)$ exists and is unique.
        Let us denote it by $z(\cdot) \coloneqq z(\cdot \mid t, w(\cdot), \nu(\cdot))$.
        Similarly, let us consider the unique solution $Z(\cdot) \coloneqq Z(\cdot \mid t, w(\cdot), \nu(\cdot)) \in \Cont^{1 - \alpha}((0, 1], \mathbb{R}^{n \times n})$ of equation \eqref{z^2_integral_equation}.

        \begin{lemma} \label{lemma_y_properties}
            Let assumptions $(f.1)$, $(f.3)$, and $(f.4)$ hold and let $N \geq c_f$, $M \geq 0$, and $\eta \in (0, T)$.
            Then, there is $\lambda_y \geq 0$ such that, for any $(t, w(\cdot)) \in G_N^\eta$, any $f \in B(M)$, any $\nu(\cdot) \in \mathcal{N}$, any $\delta \in (0, \eta / 2]$, and any $\vartheta \in (0, 1]$,
            \begin{equation} \label{lemma_y_properties_lambda_y}
                \biggl\| \frac{y(\vartheta \mid t + \delta, x^{(f)}_{t + \delta} (\cdot \mid t, w(\cdot)), \nu(\cdot))
                - y(\vartheta \mid t, w(\cdot), \nu(\cdot))}{\delta} \biggr\|
                \leq \frac{\lambda_y}{\vartheta^{1 - \alpha}},
            \end{equation}
            where $x^{(f)}_{t + \delta}(\cdot \mid t, w(\cdot))$ is the restriction of $x^{(f)}(\cdot \mid t, w(\cdot))$ {\rm(}see \eqref{x^f}{\rm)} to $[0, t + \delta]$ {\rm(}see \eqref{x_t}{\rm)} and $y(\cdot \mid t + \delta, x^{(f)}_{t + \delta} (\cdot \mid t, w(\cdot)), \nu(\cdot))$ is the solution of the integral equation \eqref{y_integral_equation}.
            Moreover,
            \begin{align*}
                & \biggl\| \frac{y(1 \mid t + \delta, x^{(f)}_{t + \delta} (\cdot \mid t, w(\cdot)), \nu(\cdot))
                - y(1 \mid t, w(\cdot), \nu(\cdot))}{\delta} \\
                & \quad - z(1 \mid t, w(\cdot), \nu(\cdot)) - Z(1 \mid t, w(\cdot), \nu(\cdot)) f \biggr\|
                \to 0
                \quad \text{as } \delta \to 0^+
            \end{align*}
            uniformly in $(t, w(\cdot)) \in G_N^\eta$, $f \in B(M)$, and $\nu(\cdot) \in \mathcal{N}$.
        \end{lemma}
        {\it Proof}
            We follow the lines of the proofs of \cite[Lemmas 2 and 4--6]{Gomoyunov_2022_FCAA}.
            Therefore, in some places, we omit details and give only the corresponding references.

            1.
                Considering the number $R_{N + M} \geq 0$, determined in accordance with \eqref{R_N}, and using property $(f^\ast.4)$, let us choose $\lambda_f^\ast \geq 0$ such that
                \begin{equation*}
                    \|f^\ast(\tau, x, \mu) - f^\ast(\tau^\prime, x^\prime, \mu)\|
                    \leq \lambda_f^\ast (|\tau - \tau^\prime| + \|x - x^\prime\|)
                \end{equation*}
                for all $\tau$, $\tau^\prime \in [0, T]$, all $x$, $x^\prime \in B(R_{N + M})$, and all $\mu \in \rpm(P)$.
                Let us put
                \begin{equation*}
                    \lambda_a
                    \coloneqq \frac{R_N + 2^{1 - \alpha} M}{\mathrm{\Gamma}(\alpha) \eta^{1 - \alpha}},
                    \quad \lambda_y^\ast
                    \coloneqq \mathrm{\Gamma}(\alpha) \lambda_a + \frac{2^{1 - \alpha} c_f (1 + R_{N + M})}{\eta^{1 - \alpha}}
                    + \frac{T^\alpha \lambda_f^\ast}{\alpha},
                \end{equation*}
                and $\lambda_y \coloneqq \lambda_y^\ast \mathrm{E}_{\alpha, \alpha} (T^\alpha \lambda_f^\ast)$, where $\mathrm{E}_{\alpha, \alpha}$ is the two-parametric Mittag-Leffler function (see, e.g., \cite[Chapter 4]{Gorenflo_Kilbas_Mainardi_Rogosin_2014}).
                Let us fix $(t, w(\cdot)) \in G_N^\eta$, $f \in B(M)$, and $\nu(\cdot) \in \mathcal{N}$.
                For any $\delta \in [0, \eta / 2]$ and any $\vartheta \in [0, 1]$, let us denote (see \eqref{a})
                \begin{align}
                    a^{(\delta)}(\vartheta)
                    & \coloneqq a(t + \delta + \vartheta (T - t - \delta) \mid t + \delta, x^{(f)}_{t + \delta} (\cdot \mid t, w(\cdot))),
                    \nonumber \\
                    y^{(\delta)}(\vartheta)
                    & \coloneqq y (\vartheta \mid t + \delta, x^{(f)}_{t + \delta} (\cdot \mid t, w(\cdot)), \nu(\cdot)).
                    \label{a^delta}
                \end{align}
                Let $\delta \in [0, \eta / 2]$ be fixed.
                Arguing as in the proof of \cite[Lemma 2]{Gomoyunov_2022_FCAA}, we derive
                \begin{align}
                    \|a^{(\delta)}(\vartheta) - a^{(0)}(\vartheta)\|
                    & \leq \frac{R_N (1 - \vartheta) \delta}{\mathrm{\Gamma}(\alpha) \vartheta^{1 - \alpha} (T - t)^{1 - \alpha}}
                    + \frac{M \delta}{\mathrm{\Gamma}(\alpha) \vartheta^{1 - \alpha} (T - t - \delta)^{1 - \alpha}}
                    \nonumber \\
                    & \leq \frac{\lambda_a \delta}{\vartheta^{1 - \alpha}}
                    \quad \forall \vartheta \in (0, 1].
                    \label{proposition_a_LD_inequality}
                \end{align}
                Further, based on \eqref{y_integral_equation}, \eqref{corollary_y_uniform_continuity_R_N+M}, and \eqref{proposition_a_LD_inequality}, similarly to the proof of \cite[Lemma 5, item (iii)]{Gomoyunov_2022_FCAA}, we obtain
                \begin{equation*}
                    \|y^{(\delta)}(\vartheta) - y^{(0)}(\vartheta)\|
                    \leq \frac{\lambda_y^\ast \delta}{\mathrm{\Gamma}(\alpha) \vartheta^{1 - \alpha}}
                    + \frac{T^\alpha \lambda_f^\ast}{\mathrm{\Gamma}(\alpha)} \int_{0}^{\vartheta}
                    \frac{\|y^{(\delta)}(\zeta) - y^{(0)}(\zeta)\|}{(\vartheta - \zeta)^{1 - \alpha}} \, \rd \zeta
                \end{equation*}
                for all $\vartheta \in (0, 1]$.
                Hence, applying a Gronwall type inequality (see, e.g., \cite[Corollary 1]{Gomoyunov_2022_FCAA}), we get
                \begin{equation*}
                    \|y^{(\delta)}(\vartheta) - y^{(0)}(\vartheta)\|
                    \leq \frac{\lambda_y^\ast \mathrm{E}_{\alpha, \alpha} (T^\alpha \lambda_f^\ast) \delta}{\vartheta^{1 - \alpha}}
                    = \frac{\lambda_y \delta}{\vartheta^{1 - \alpha}}
                    \quad \forall \vartheta \in (0, 1],
                \end{equation*}
                which gives \eqref{lemma_y_properties_lambda_y}.

            2.
                Let us prove the second part of the lemma.
                Let us fix $(t, w(\cdot)) \in G_N^\eta$, $f \in B(M)$, and $\nu(\cdot) \in \mathcal{N}$.
                Let us consider the functions $q(\cdot)$ and $Q(\cdot)$ given by \eqref{r_R} and the solutions $z(\cdot) \coloneqq z(\cdot \mid t, w(\cdot), \nu(\cdot))$ and $Z(\cdot) \coloneqq Z(\cdot \mid t, w(\cdot), \nu(\cdot))$ of the integral equations \eqref{z^1_integral_equation} and \eqref{z^2_integral_equation}, respectively.
                For any $\delta \in (0, \eta / 2]$ and any $\vartheta \in (0, 1]$, using the notation introduced in \eqref{a^delta}, let us put
                \begin{align*}
                    r_a^{(\delta)}(\vartheta)
                    & \coloneqq \biggl\| \frac{a^{(\delta)}(\vartheta) - a^{(0)}(\vartheta)}{\delta} - q(\vartheta) - Q(\vartheta) f \biggr\|, \\
                    r_y^{(\delta)}(\vartheta)
                    & \coloneqq \biggl\| \frac{y^{(\delta)}(\vartheta) - y^{(0)}(\vartheta)}{\delta} - z(\vartheta) - Z(\vartheta) f \biggr\|.
                \end{align*}
                Let $\delta \in (0, \eta / 2]$ and $\vartheta \in (0, 1]$ be fixed.
                Arguing as in the proof of \cite[Lemma 4, equality (5.5)]{Gomoyunov_2022_FCAA}, we obtain
                \begin{equation*}
                    r_a^{(\delta)}(\vartheta)
                    \leq \frac{R_N + 2 M}{\mathrm{\Gamma}(\alpha)}
                    \biggl( \frac{1}{\vartheta^{1 - \alpha} (T - t - \delta)^{1 - \alpha}}
                    - \frac{1}{(\delta + \vartheta (T - t - \delta))^{1 - \alpha}} \biggr).
                \end{equation*}
                Since the function in the parentheses is increasing in $t \in [0, T - \eta]$, we have
                \begin{equation} \label{proof_lemma_a_2}
                    r_a^{(\delta)}(\vartheta)
                    \leq \frac{R_N + 2 M}{\mathrm{\Gamma}(\alpha)} \biggl( \frac{1}{\vartheta^{1 - \alpha} (\eta - \delta)^{1 - \alpha}}
                    - \frac{1}{(\delta + \vartheta (\eta - \delta))^{1 - \alpha}} \biggr)
                \end{equation}
                and, in particular,
                \begin{equation} \label{proof_lemma_a_2_1}
                    r_a^{(\delta)}(1)
                    \leq \frac{R_N + 2 M}{\mathrm{\Gamma}(\alpha)} \biggl( \frac{1}{(\eta - \delta)^{1 - \alpha}}
                    - \frac{1}{\eta^{1 - \alpha}} \biggr).
                \end{equation}
                Based on \eqref{proof_lemma_a_2}, similarly to the proof of \cite[Lemma 4, equality (5.6)]{Gomoyunov_2022_FCAA}, we derive
                \begin{align}
                    & \int_{0}^{1}
                    \frac{r_a^{(\delta)}(\zeta)}{(1 - \zeta)^{1 - \alpha}} \, \rd \zeta
                    \nonumber \\
                    & \ \leq \frac{R_N + 2 M}{\mathrm{\Gamma}(\alpha) (\eta - \delta)^{1 - \alpha}}
                    \biggl( \int_{0}^{1} \frac{\rd \zeta}{(1 - \zeta)^{1 - \alpha} \zeta^{1 - \alpha}}
                    - \int_{0}^{1} \frac{\rd \zeta}
                    {(1 - \zeta)^{1 - \alpha} \bigl( \frac{\delta}{\eta - \delta} + \zeta \bigr)^{1 - \alpha}} \biggr)
                    \nonumber \\
                    & \ \leq \frac{R_N + 2 M}{\mathrm{\Gamma}(\alpha) (\eta - \delta)^{1 - \alpha}}
                    \biggl( \mathrm{B}(\alpha, \alpha) \biggl( 1
                    - \biggl( 1 + \frac{\delta}{\eta - \delta} \biggr)^{2 \alpha - 1} \biggr)
                    + \frac{1}{\alpha} \biggl( \frac{\delta}{\eta - \delta} \biggr)^\alpha \biggr),
                    \label{proof_lemma_a_3}
                \end{align}
                where $\mathrm{B}$ is the beta-function.
                Let us note that the functions from the right-hand sides of estimates \eqref{proof_lemma_a_2_1} and \eqref{proof_lemma_a_3} do not depend on a particular choice of $(t, w(\cdot))$ and $f$ and tend to zero as $\delta \to 0^+$.
                Hence, we conclude that
                \begin{equation} \label{a_limit_relations}
                    r_a^{(\delta)}(1)
                    \to 0,
                    \quad \int_{0}^{1}
                    \frac{r_a^{(\delta)}(\zeta)}{(1 - \zeta)^{1 - \alpha}} \, \rd \zeta
                    \to 0
                    \quad \text{as } \delta \to 0^+
                \end{equation}
                uniformly in $(t, w(\cdot)) \in G_N^\eta$ and $f \in B(M)$.

            3.
                Let us consider the set $\Omega \coloneqq [0, T - \eta / 2] \times [0, 1] \times B(R_{N + M}) \times \rpm(P)$ and the function
                \begin{equation*}
                    \omega^\ast (\tau, \vartheta, x, \mu)
                    \coloneqq (T - \tau)^\alpha f^\ast(\tau + \vartheta (T - \tau), x, \mu)
                    \quad \forall (\tau, \vartheta, x, \mu) \in \Omega.
                \end{equation*}
                By properties $(f^\ast.1)$ and $(f^\ast.4)$, the function $\omega^\ast$ is continuous together with its partial derivatives (in this connection, see \eqref{f_ast}, \eqref{f_ast_partial_derivatives}, \eqref{A_b}, and \eqref{A_b^ast})
                \begin{align*}
                    \frac{\partial \omega^\ast}{\partial \tau} (\tau, \vartheta, x, \mu)
                    & = (1 - \vartheta) (T - \tau)^\alpha \frac{\partial f^\ast}{\partial \tau} (\tau + \vartheta (T - \tau), x, \mu) \\
                    & \quad - \frac{\alpha f^\ast(\tau + \vartheta (T - \tau), x, \mu)}{(T - \tau)^{1 - \alpha}}, \\
                    \frac{\partial \omega^\ast}{\partial x} (\tau, \vartheta, x, \mu)
                    & = (T - \tau)^\alpha \frac{\partial f^\ast}{\partial x} (\tau + \vartheta (T - \tau), x, \mu)
                    \quad \forall (\tau, \vartheta, x, \mu) \in \Omega.
                \end{align*}
                Since the set $\Omega$ is compact, there exists $R \geq 0$ such that
                \begin{equation*}
                    \biggl\| \frac{\partial \omega^\ast}{\partial x} (\tau, \vartheta, x, \mu) \biggr\|
                    \leq R
                    \quad \forall (\tau, \vartheta, x, \mu) \in \Omega.
                \end{equation*}

            4.
                Now, let $\varepsilon > 0$ be given.
                Due to the uniform convergence \eqref{a_limit_relations}, there exists $\delta_1 \in (0, \eta / 2]$ such that, for any $(t, w(\cdot)) \in G_N^\eta$, any $f \in B(M)$, and any $\delta \in (0, \delta_1]$,
                \begin{equation*}
                    r_a^{(\delta)}(1)
                    + R \mathrm{E}_{\alpha, \alpha} (R) \int_{0}^{1} \frac{r_a^{(\delta)}(\zeta)}{(1 - \zeta)^{1 - \alpha}} \, \rd \zeta
                    \leq \frac{\varepsilon}{2}.
                \end{equation*}
                Let us choose $\kappa > 0$ from the condition
                \begin{equation*}
                    \kappa \mathrm{E}_{\alpha, \alpha} (R) (1 + \lambda_y) \mathrm{B}(\alpha, \alpha)
                    \leq \frac{\varepsilon}{2}
                \end{equation*}
                and take $\chi > 0$ such that, for any $(\tau, \vartheta, x, \mu)$, $(\tau^\prime, \vartheta, x^\prime, \mu) \in \Omega$ with $|\tau - \tau^\prime| \leq \chi$ and $\|x - x^\prime\| \leq \chi$,
                \begin{align*}
                    \biggl\| \frac{\partial \omega^\ast}{\partial \tau} (\tau, \vartheta, x, \mu)
                    - \frac{\partial \omega^\ast}{\partial \tau} (\tau^\prime, \vartheta, x^\prime, \mu) \biggr\|
                    & \leq \kappa, \\
                    \quad \biggl\| \frac{\partial \omega^\ast}{\partial x} (\tau, \vartheta, x, \mu)
                    - \frac{\partial \omega^\ast}{\partial x} (\tau^\prime, \vartheta, x^\prime, \mu) \biggr\|
                    & \leq \kappa.
                \end{align*}
                By Corollary \ref{corollary_y_uniform_continuity}, there exists $\delta_2 \in (0, \eta / 2]$ such that, for any $(t, w(\cdot)) \in G_N^\eta$, any $f \in B(M)$, any $\nu(\cdot) \in \mathcal{N}$, and any $\delta \in (0, \delta_2]$,
                \begin{equation*}
                    \|y^{(\delta)}(\cdot) - y^{(0)}(\cdot)\|_{[0, 1]}
                    \leq \chi.
                \end{equation*}
                Let us put $\delta_\ast \coloneqq \min\{\delta_1, \delta_2, \chi\} > 0$.

                Let $(t, w(\cdot)) \in G_N^\eta$, $f \in B(M)$, $\nu(\cdot) \in \mathcal{N}$, and $\delta \in (0, \delta_\ast]$ be fixed.
                Then, arguing as in the proof of \cite[Lemma 6]{Gomoyunov_2022_FCAA}, we derive
                \begin{equation*}
                    r_y^{(\delta)}(\vartheta)
                    \leq r_a^{(\delta)}(\vartheta)
                    + \frac{\kappa (1 + \lambda_y) \mathrm{B}(\alpha, \alpha)}{\mathrm{\Gamma}(\alpha) \vartheta^{1 - \alpha}}
                    + \frac{R}{\mathrm{\Gamma}(\alpha)} \int_{0}^{\vartheta} \frac{r_y^{(\delta)}(\zeta)}{(\vartheta - \zeta)^{1 - \alpha}} \, \rd \zeta
                    \quad \forall \vartheta \in (0, 1].
                \end{equation*}
                Hence, noting that $r_a^{(\delta)}(\cdot)$, $r_y^{(\delta)}(\cdot) \in \Cont^{1 - \alpha}((0, 1], \mathbb{R})$ and applying a Gronwall type inequality (see, e.g., \cite[Corollary 1]{Gomoyunov_2022_FCAA}), we obtain
                \begin{align*}
                    r_y^{(\delta)}(1)
                    & \leq \kappa \mathrm{E}_{\alpha, \alpha}(R) (1 + \lambda_y) \mathrm{B}(\alpha, \alpha) + r_a^{(\delta)}(1)
                    + R \mathrm{E}_{\alpha, \alpha}(R) \int_{0}^{1} \frac{r_a^{(\delta)}(\zeta)}{(1 - \zeta)^{1 - \alpha}} \, \rd \zeta \\
                    & \leq \varepsilon,
                \end{align*}
                which completes the proof of the lemma.
        \qed

        Let us prove a continuity dependence property of solutions $z(\cdot \mid t, w(\cdot), \nu(\cdot))$ and $Z(\cdot \mid t, w(\cdot), \nu(\cdot))$ of the integral equations \eqref{z^1_integral_equation} and \eqref{z^2_integral_equation} with respect to the parameter $\nu(\cdot) \in \mathcal{N}$.
        \begin{lemma} \label{lemma_derivatives_of_y_are_continuous}
            Let assumptions $(f.1)$, $(f.3)$, and $(f.4)$ hold.
            Then, for every $(t, w(\cdot)) \in G^0$, the following mappings are continuous:
            \begin{align}
                \mathcal{N} \ni \nu(\cdot)
                & \mapsto z(1 \mid t, w(\cdot), \nu(\cdot)) \in \mathbb{R}^n,
                \nonumber \\
                \mathcal{N} \ni \nu(\cdot)
                & \mapsto Z(1 \mid t, w(\cdot), \nu(\cdot)) \in \mathbb{R}^{n \times n}.
                \label{mappings_derivatives}
            \end{align}
        \end{lemma}
        {\it Proof}
            We prove continuity of the first mapping from \eqref{mappings_derivatives} only, since continuity of the second one can be verified in a similar way.
            Let $\{\nu_i(\cdot)\}_{i \in \mathbb{N}_0} \subset \mathcal{N}$ be such that $\nu_i(\cdot) \to \nu_0(\cdot)$ as $i \to \infty$ and let $z_i(\cdot) \coloneqq z(\cdot \mid t, w(\cdot), \nu_i(\cdot))$ for all $i \in \mathbb{N}_0$.
            It is required to show that $z_i(1) \to z_0(1)$ as $i \to \infty$.

            For every $i \in \mathbb{N}_0$, let us consider the solution $y_i(\cdot) \coloneqq y(\cdot \mid t, w(\cdot), \nu_i(\cdot))$ of the integral equation \eqref{y_integral_equation} and introduce the corresponding functions $A_i$, $b_i$ and $A_i^\ast$, $b_i^\ast$ according to \eqref{A_b} and \eqref{A_b^ast}, respectively.
            Due to assumptions $(f.1)$ and $(f.4)$ and since $y_i(\cdot) \to y_0(\cdot)$ as $i \to \infty$ by Lemma \ref{lemma_y_continuous}, we obtain that $\|A_i(\vartheta, u) - A_0(\vartheta, u)\| \to 0$ and $\|b_i(\vartheta, u) - b_0(\vartheta, u)\| \to 0$ as $i \to \infty$ uniformly in $\vartheta \in [0, 1]$ and $u \in P$, and, hence,
            \begin{equation} \label{uniform_convergence}
                \| A_i^\ast(\vartheta, \mu) - A_0^\ast(\vartheta, \mu)\| \to
                0,
                \quad \|b_i^\ast(\vartheta, \mu) - b_0^\ast(\vartheta, \mu)\| \to
                0 \quad \text{as } i \to \infty
            \end{equation}
            uniformly in $\vartheta \in [0, 1]$ and $\mu \in \rpm(P)$.
            In particular, there is $R \geq 0$ such that
            \begin{equation*}
                \|A_i^\ast(\vartheta, \mu)\|
                \leq R
                \quad \forall \vartheta \in [0, 1] \quad \forall \mu \in \rpm(P) \quad \forall i \in \mathbb{N}_0.
            \end{equation*}

            Based on \eqref{z^1_integral_equation}, for any $i \in \mathbb{N}_0$ and any $\vartheta \in (0, 1]$, we derive
            \begin{align*}
                \|z_i(\vartheta) - z_0(\vartheta)\|
                & = \biggl\| \frac{1}{\mathrm{\Gamma}(\alpha)} \int_{0}^{\vartheta}
                \frac{A_i^\ast(\zeta, \nu_i(\zeta)) z_i(\zeta) + b_i^\ast(\zeta, \nu_i(\zeta))}
                {(\vartheta - \zeta)^{1 - \alpha}} \, \rd \zeta \\
                & \quad - \frac{1}{\mathrm{\Gamma}(\alpha)} \int_{0}^{\vartheta}
                \frac{A_0^\ast(\zeta, \nu_0(\zeta)) z_0(\zeta) + b_0^\ast(\zeta, \nu_0(\zeta))}
                {(\vartheta - \zeta)^{1 - \alpha}} \, \rd \zeta \biggr\| \\
                & \leq \frac{R}{\mathrm{\Gamma}(\alpha)}
                \int_{0}^{\vartheta} \frac{\|z_i(\zeta) -  z_0(\zeta)\|}{(\vartheta - \zeta)^{1 - \alpha}} \, \rd \zeta
                + r^{(1)}_i(\vartheta) + r^{(2)}_i(\vartheta),
            \end{align*}
            where we denote
            \begin{align*}
                r^{(1)}_i(\vartheta)
                & \coloneqq \frac{1}{\mathrm{\Gamma}(\alpha)} \int_{0}^{\vartheta}
                \frac{\|A_i^\ast(\zeta, \nu_i(\zeta)) - A_0^\ast(\zeta, \nu_i(\zeta))\| \|z_0(\zeta)\|}
                {(\vartheta - \zeta)^{1 - \alpha}} \, \rd \zeta \\
                & \quad + \frac{1}{\mathrm{\Gamma}(\alpha)} \int_{0}^{\vartheta}
                \frac{\|b_i^\ast(\zeta, \nu_i(\zeta)) - b_0^\ast(\zeta, \nu_i(\zeta))\|}{(\vartheta - \zeta)^{1 - \alpha}} \, \rd \zeta, \\
                r^{(2)}_i(\vartheta)
                & \coloneqq \biggl\| \frac{1}{\mathrm{\Gamma}(\alpha)}
                \int_{0}^{\vartheta} \frac{A_0^\ast(\zeta, \nu_i(\zeta)) z_0(\zeta) + b_0^\ast (\zeta, \nu_i(\zeta))}
                {(\vartheta - \zeta)^{1 - \alpha}} \rd \zeta \\
                & \quad - \frac{1}{\mathrm{\Gamma}(\alpha)}
                \int_{0}^{\vartheta} \frac{A_0^\ast(\zeta, \nu_0(\zeta)) z_0(\zeta) + b_0^\ast (\zeta, \nu_0(\zeta))}
                {(\vartheta - \zeta)^{1 - \alpha}} \, \rd \zeta \biggr\|.
            \end{align*}

            Now, let $\varepsilon > 0$ be given.
            Let us choose $\varepsilon_\ast > 0$ from the condition
            \begin{equation*}
                \varepsilon_\ast \mathrm{E}_{\alpha, \alpha} (R) \biggl( \mathrm{B}(\alpha, \alpha) + \frac{1}{\alpha} + 1 \biggr)
                \leq \varepsilon,
            \end{equation*}
            where $\mathrm{E}_{\alpha, \alpha}$ denotes again the two-parametric Mittag-Leffler function.
            In view of the inclusion $z_0(\cdot) \in \Cont^{1 - \alpha}((0, 1], \mathbb{R}^n)$ and the uniform convergence \eqref{uniform_convergence}, there exists $i_1 \in \mathbb{N}$ such that, for any $i \in \mathbb{N}$ with $i \geq i_1$ and any $\vartheta \in (0, 1]$,
            \begin{align*}
                r_i^{(1)}(\vartheta)
                & \leq \frac{\varepsilon_\ast}{\mathrm{\Gamma}(\alpha)} \int_{0}^{\vartheta} \frac{\rd \zeta}{(\vartheta - \zeta)^{1 - \alpha} \zeta^{1 - \alpha}}
                + \frac{\varepsilon_\ast}{\mathrm{\Gamma}(\alpha)} \int_{0}^{\vartheta} \frac{\rd \zeta}{(\vartheta - \zeta)^{1 - \alpha}} \\
                & \leq \frac{\varepsilon_\ast}{\mathrm{\Gamma}(\alpha) \vartheta^{1 - \alpha}} \biggl( \mathrm{B}(\alpha, \alpha) + \frac{1}{\alpha} \biggr).
            \end{align*}
            Further, for every $i \in \mathbb{N}_0$, let us consider the auxiliary function
            \begin{equation*}
                h_i(\vartheta)
                \coloneqq \frac{\vartheta^{1 - \alpha}}{\mathrm{\Gamma}(\alpha)}
                \int_{0}^{\vartheta} \frac{A_0^\ast(\zeta, \nu_i(\zeta)) z_0(\zeta) + b_0^\ast (\zeta, \nu_i(\zeta))}
                {(\vartheta - \zeta)^{1 - \alpha}} \rd \zeta
                \quad \forall \vartheta \in [0, 1].
            \end{equation*}
            Let us fix $\vartheta \in (0, 1]$ and introduce the function $\mathfrak{b}^{[\vartheta]} \colon [0, 1] \times P \to \mathbb{R}^n$ by
            \begin{equation*}
                \mathfrak{b}^{[\vartheta]} (\zeta, u)
                \coloneqq
                \begin{cases}
                    \displaystyle
                    \frac{\vartheta^{1 - \alpha} (A_0(\zeta, u) z_0(\zeta) + b_0(\zeta, u))}{\mathrm{\Gamma}(\alpha)(\vartheta - \zeta)^{1 - \alpha}},
                    & \mbox{if } \zeta \in (0, \vartheta), \\
                    0,
                    & \mbox{if } \zeta \in \{0\} \cup [\vartheta, 1],
                \end{cases}
            \end{equation*}
            for all $u \in P$.
            Note that the function $\mathfrak{b}^{[\vartheta]}$ is a Carath\'{e}odory function due to continuity of the functions $A_0$ and $b_0$ and the inclusion $z_0(\cdot) \in \Cont((0, 1], \mathbb{R}^n)$.
            Therefore, by definition \eqref{convergence_M} of convergence in $\mathcal{N}$, we have
            \begin{equation*}
                h_i(\vartheta)
                = \int_{0}^{1} \int_{P} \mathfrak{b}^{[\vartheta]} (\zeta, u) \nu_i(\zeta, \rd u) \, \rd \zeta
                \to \int_{0}^{1} \int_{P} \mathfrak{b}^{[\vartheta]} (\zeta, u) \nu_0(\zeta, \rd u) \, \rd \zeta
                = h_0(\vartheta)
            \end{equation*}
            as $i \to \infty$.
            Consequently, and since $h_i(0) = 0$ for all $i \in \mathbb{N}_0$, we conclude that $h_i(\vartheta) \to h_0(\vartheta)$ as $i \to \infty$ for all $\vartheta \in [0, 1]$.
            In addition, by continuity of the functions $A_0^\ast$ and $b_0^\ast$ and the inclusion $z_0(\cdot) \in \Cont((0, 1], \mathbb{R}^n)$, there exists $\varkappa \geq 0$ such that $\zeta^{1 - \alpha} \| A_0^\ast (\zeta, \mu) z_0(\zeta) + b_0^\ast(\zeta, \mu)\| \leq \varkappa$ for all $\zeta \in (0, 1]$ and all $\mu \in \rpm(P)$, and, then, it follows from \cite[Corollary 2.1]{Gomoyunov_2019_FCAA_2} that
            \begin{equation*}
                \|h_i(\vartheta) - h_i(\vartheta^\prime)\|
                \leq \frac{2 \varkappa}{\mathrm{\Gamma}(\alpha + 1)} \biggl(1 + \frac{\sin(\alpha \pi)}{\alpha \pi} \biggr) |\vartheta - \vartheta^\prime|^\alpha
            \end{equation*}
            for all $i \in \mathbb{N}_0$ and all $\vartheta$, $\vartheta^\prime \in [0, 1]$.
            Hence, the functions $h_i(\cdot)$ for all $i \in \mathbb{N}_0$ are equicontinuous, which together with the pointwise convergence established above implies that $\{h_i(\cdot)\}_{i \in \mathbb{N}}$ converges to $h_0(\cdot)$ as $i \to \infty$ uniformly on $[0, 1]$ (see, e.g., \cite[Theorem I.5.3]{Warga_1972}).
            Thus, there exists $i_2 \in \mathbb{N}$ such that, for any $i \in \mathbb{N}$ with $i \geq i_2$ and any $\vartheta \in (0, 1]$,
            \begin{equation*}
                r_i^{(2)}(\vartheta)
                = \frac{\|h_i(\vartheta) - h_0(\vartheta)\|}{\vartheta^{1 - \alpha}}
                \leq \frac{\varepsilon_\ast}{\mathrm{\Gamma}(\alpha) \vartheta^{1 - \alpha}}.
            \end{equation*}

            As a result, for any $i \in \mathbb{N}$ with $i \geq \max\{i_1, i_2\}$, we get
            \begin{equation*}
                \|z_i(\vartheta) - z_0(\vartheta)\|
                \leq \frac{\varepsilon_\ast}{\mathrm{\Gamma}(\alpha) \vartheta^{1 - \alpha}} \biggl( \mathrm{B}(\alpha, \alpha) + \frac{1}{\alpha} + 1 \biggr)
                + \frac{R}{\mathrm{\Gamma}(\alpha)} \int_{0}^{\vartheta} \frac{\|z_i(\zeta) - z_0(\zeta)\|}{(\vartheta - \zeta)^{1 - \alpha}} \, \rd \zeta
            \end{equation*}
            for all $\vartheta \in (0, 1]$, wherefrom, recalling that $z_i(\cdot)$, $z_0(\cdot) \in \Cont((0, 1], \mathbb{R}^n)$ and applying a Gronwall type inequality (see, e.g., \cite[Corollary 1]{Gomoyunov_2022_FCAA}), we derive
            \begin{equation*}
                \|z_i(1) - z_0(1)\|
                \leq \varepsilon_\ast \mathrm{E}_{\alpha, \alpha} (R) \biggl( \mathrm{B}(\alpha, \alpha) + \frac{1}{\alpha} + 1 \biggr)
                \leq \varepsilon,
            \end{equation*}
            which completes the proof.
        \qed

    \subsection{Directional Differentiability of Functional $\psi$}

        From now on, let us assume that the function $\sigma$ from the cost functional \eqref{cost_functional} satisfies the following condition:
        \begin{description}
            \item[$(\sigma.2)$]
                The partial derivatives $\partial \sigma / \partial x \colon \mathbb{R}^n \to \mathbb{R}^n$ exist and are continuous.
        \end{description}
        Let us note that condition $(\sigma.2)$ implies condition $(\sigma.1)$.

        Let us consider the functional $\psi$ given by \eqref{psi} and put
        \begin{align}
            \partial^\alpha_t \psi(t, w(\cdot), \nu(\cdot))
            & \coloneqq \biggl\langle \frac{\partial \sigma}{\partial x} \bigl( y(1 \mid t, w(\cdot), \nu(\cdot)) \bigr),
            z(1 \mid t, w(\cdot), \nu(\cdot)) \biggr\rangle,
            \nonumber \\
            \nabla^\alpha \psi(t, w(\cdot), \nu(\cdot))
            & \coloneqq Z(1 \mid t, w(\cdot), \nu(\cdot))^\top
            \frac{\partial \sigma}{\partial x} \bigl( y(1 \mid t, w(\cdot), \nu(\cdot)) \bigr)
            \label{psi_formulas_for_derivatives}
        \end{align}
        for all $(t, w(\cdot)) \in G^0$ and all $\nu(\cdot) \in \mathcal{N}$.
        In the above, $y(\cdot \mid t, w(\cdot), \nu(\cdot))$ is the solution of the integral equation \eqref{y_integral_equation}, $z(\cdot \mid t, w(\cdot), \nu(\cdot))$ and $Z(\cdot \mid t, w(\cdot), \nu(\cdot))$ are the solutions of the integral equations \eqref{z^1_integral_equation} and \eqref{z^2_integral_equation}, respectively, and the superscript $^\top$ denotes transposition.

        \begin{lemma} \label{lemma_psi_is_differentiable}
            Let assumptions $(f.1)$, $(f.3)$, $(f.4)$, and $(\sigma.2)$ hold.
            Then, for any $N \geq c_f$, any $M \geq 0$, and any $\eta \in (0, T)$,
            \begin{align*}
                & \biggl| \frac{\psi(t + \delta, x^{(f)}_{t + \delta}(\cdot \mid t, w(\cdot)), \nu(\cdot)) - \psi(t, w(\cdot), \nu(\cdot))}{\delta} \\
                & \quad - \partial^\alpha_t \psi(t, w(\cdot), \nu(\cdot)) - \langle \nabla^\alpha \psi(t, w(\cdot), \nu(\cdot)), f \rangle \biggr|
                \to 0
                \quad \text{as } \delta \to 0^+
            \end{align*}
            uniformly in $(t, w(\cdot)) \in G_N^\eta$, $f \in B(M)$, and $\nu(\cdot) \in \mathcal{N}$, where $x^{(f)}_{t + \delta}(\cdot \mid t, w(\cdot))$ is the restriction of $x^{(f)}(\cdot \mid t, w(\cdot))$ {\rm(}see \eqref{x^f}{\rm)} to $[0, t + \delta]$ {\rm(}see \eqref{x_t}{\rm)}.
        \end{lemma}
        {\it Proof}
            Let us take $\lambda_y \geq 0$ from Lemma \ref{lemma_y_properties}.
            Let us consider the number $R_{N + M} \geq 0$, determined in accordance with \eqref{R_N}, and, based on assumption $(\sigma.2)$, choose $R \geq 0$ such that
            \begin{equation*}
                \biggl\| \frac{\partial \sigma}{\partial x}(x) \biggr\|
                \leq R
                \quad \forall x \in B(R_{N + M}).
            \end{equation*}

            Let $\varepsilon > 0$ be given.
            Let us choose $\varepsilon_\ast > 0$ from the condition $\varepsilon_\ast (\lambda_y + R) \leq \varepsilon$.
            By Lemma \ref{lemma_y_properties}, there exists $\delta_1 \in (0, \eta / 2]$ such that, for any $(t, w(\cdot)) \in G_N^\eta$, any $f \in B(M)$, any $\nu(\cdot) \in \mathcal{N}$, and any $\delta \in (0, \delta_1]$,
            \begin{equation*}
                \bigg\| \frac{y^{(\delta)}(1) - y^{(0)}(1)}{\delta}
                - z(1) - Z(1) f \bigg\|
                \leq \varepsilon_\ast,
            \end{equation*}
            where $z(1) \coloneqq z(1 \mid t, w(\cdot), \nu(\cdot))$ and $Z(1) \coloneqq Z(1 \mid t, w(\cdot), \nu(\cdot))$ and the notation introduced in \eqref{a^delta} is used.
            Using assumption $(\sigma.2)$, let us take $\kappa > 0$ such that, for any $x$, $x^\prime \in B(R_{N + M})$ with $\|x - x^\prime\| \leq \kappa$,
            \begin{equation*}
                \biggl\| \frac{\partial \sigma}{\partial x}(x) - \frac{\partial \sigma}{\partial x}(x^\prime) \biggr\|
                \leq \varepsilon_\ast.
            \end{equation*}
            By Corollary \ref{corollary_y_uniform_continuity}, there exists $\delta_2 \in (0, \eta / 2]$ such that, for any $(t, w(\cdot)) \in G_N^\eta$, any $f \in B(M)$, any $\nu(\cdot) \in \mathcal{N}$, and any $\delta \in (0, \delta_2]$,
            \begin{equation*}
                \|y^{(\delta)}(1) - y^{(0)}(1)\|
                \leq \kappa.
            \end{equation*}
            Let us put $\delta_\ast \coloneqq \min\{\delta_1, \delta_2\} > 0$.

            Let $(t, w(\cdot)) \in G_N^\eta$, $f \in B(M)$, $\nu(\cdot) \in \mathcal{N}$, and $\delta \in (0, \delta_\ast]$ be fixed.
            We have
            \begin{equation*}
                \psi(t + \delta, x^{(f)}_{t + \delta}(\cdot \mid t, w(\cdot)), \nu(\cdot))
                - \psi(t, w(\cdot), \nu(\cdot))
                = \sigma(y^{(\delta)}(1)) - \sigma(y^{(0)}(1)).
            \end{equation*}
            By the mean value theorem, on the segment connecting the points $y^{(\delta)}(1)$ and $y^{(0)}(1)$, there exists a point $y^\prime$ for which
            \begin{equation*}
                \sigma(y^{(\delta)}(1)) - \sigma(y^{(0)}(1))
                = \biggl\langle \frac{\partial \sigma}{\partial x} (y^\prime), y^{(\delta)}(1) - y^{(0)}(1) \biggr\rangle.
            \end{equation*}
            Since $y^{(\delta)}(1)$, $y^{(0)}(1) \in B(R_{N + M})$ by Corollary \ref{corollary_y_uniform_continuity}, we get $y^\prime \in B(R_{N + M})$.
            Moreover, we derive $\|y^{(0)}(1) - y^\prime\| \leq \|y^{(\delta)}(1) - y^{(0)}(1)\| \leq \kappa$ and, consequently,
            \begin{equation*}
                \biggl\| \frac{\partial \sigma}{\partial x} (y^{(0)}(1)) - \frac{\partial \sigma}{\partial x} (y^\prime) \biggr\|
                \leq \varepsilon_\ast.
            \end{equation*}
            Hence, taking \eqref{psi_formulas_for_derivatives} into account, we get
            \begin{align*}
                & \biggl| \frac{\sigma(y^{(\delta)}(1)) - \sigma(y^{(0)}(1))}{\delta}
                - \partial^\alpha_t \psi(t, w(\cdot), \nu(\cdot)) - \langle \nabla^\alpha \psi(t, w(\cdot), \nu(\cdot)), f \rangle \biggr| \\
                & \quad = \biggl| \biggl\langle \frac{\partial \sigma}{\partial x} (y^\prime), \frac{y^{(\delta)}(1) - y^{(0)}(1)}{\delta} \biggr\rangle
                - \biggl\langle \frac{\partial \sigma}{\partial x} (y^{(0)}(1)), z(1) + Z(1)f \biggr\rangle \biggr| \\
                & \quad \leq \biggl\| \frac{\partial \sigma}{\partial x} (y^\prime) - \frac{\partial \sigma}{\partial x} (y^{(0)}(1)) \biggr\|
                \biggl\| \frac{y^{(\delta)}(1) - y^{(0)}(1)}{\delta} \biggr\| \\
                & \qquad + \biggl\| \frac{\partial \sigma}{\partial x} (y^{(0)}(1)) \biggr\|
                \biggl\| \frac{y^{(\delta)}(1) - y^{(0)}(1)}{\delta} - z(1) - Z(1) f \biggr\| \\
                & \quad \leq \varepsilon_\ast \lambda_y + \varepsilon_\ast R
                \leq \varepsilon.
            \end{align*}
            The proof is complete.
        \qed

        In addition, let us note that Lemmas \ref{lemma_y_continuous} and \ref{lemma_derivatives_of_y_are_continuous} immediately imply
        \begin{corollary} \label{corollary_derivatives_of_psi_are_continuous}
            Let assumptions $(f.1)$, $(f.3)$, $(f.4)$, and $(\sigma.2)$ hold.
            Then, for every $(t, w(\cdot)) \in G^0$, the following mappings are continuous:
            \begin{equation*}
                \mathcal{N} \ni \nu(\cdot) \mapsto \partial^\alpha_t \psi(t, w(\cdot), \nu(\cdot)) \in \mathbb{R},
                \quad \mathcal{N} \ni \nu(\cdot) \mapsto \nabla^\alpha \psi(t, w(\cdot), \nu(\cdot)) \in \mathbb{R}^{n}.
            \end{equation*}
        \end{corollary}

\section{Directional Differentiability of Order $\alpha$ of Value Functional}
\label{section_Directional_differentiability}

    Based on representation \eqref{value_functional_representation} of the value functional $\rho$ of the optimal control problem \eqref{system}, \eqref{cost_functional}, using compactness of the metric space $\mathcal{N}$ and the properties of the functional $\psi$ from \eqref{psi} established in Sections \ref{section_continuity} and \ref{section_differentiability} (see Corollaries \ref{corollary_psi_is_continuous} and \ref{corollary_derivatives_of_psi_are_continuous} and Lemma \ref{lemma_psi_is_differentiable}), and applying Theorem \ref{theorem_envelope_theorem}, we arrive at the following result.
    \begin{theorem} \label{theorem_directional_differentiability}
        Under assumptions $(f.1)$, $(f.3)$, $(f.4)$, and $(\sigma.2)$, the value functional $\rho$ of the optimal control problem \eqref{system}, \eqref{cost_functional} is directionally differentiable of order $\alpha$.
        In addition, for any $(t, w(\cdot)) \in G^0$ and any $f \in \mathbb{R}^n$,
        \begin{align}
            & \partial^\alpha \{ \rho(t, w(\cdot)) \mid f \}
            \nonumber \\
            & \quad = \min_{\nu(\cdot) \in \mathcal{N}^\circ(t, w(\cdot))}
            \bigl( \partial^\alpha_t \psi(t, w(\cdot), \nu(\cdot))
            + \langle \nabla^\alpha \psi(t, w(\cdot), \nu(\cdot)), f \rangle \bigr),
            \label{directional_derivative_formula}
        \end{align}
        where $\partial^\alpha_t \psi(t, w(\cdot), \nu(\cdot))$ and $\nabla^\alpha \psi(t, w(\cdot), \nu(\cdot))$ are given by \eqref{psi_formulas_for_derivatives} and
        \begin{equation} \label{N^0}
            \mathcal{N}^\circ(t, w(\cdot))
            \coloneqq \bigl\{ \nu(\cdot) \in \mathcal{N} \colon
            \psi(t, w(\cdot), \nu(\cdot)) = \rho(t, w(\cdot)) \bigr\}.
        \end{equation}
    \end{theorem}

    In particular, as a consequence of this result and Theorem \ref{theorem_non-smooth_general}, we obtain
    \begin{corollary} \label{corollary_non-convex}
        Under assumptions $(f.1)$, $(f.3)$, $(f.4)$, and $(\sigma.2)$, a functional $\varphi \colon G \to \mathbb{R}$ is the value functional of the optimal control problem \eqref{system}, \eqref{cost_functional} if and only if $\varphi$ is continuous, possesses the local Lipschitz continuity property $(L)$, is directionally differentiable of order $\alpha$, meets the boundary condition \eqref{boundary_condition}, and, for any $(t, w(\cdot)) \in G^0$, satisfies the pair of differential inequalities
        \begin{equation*}
            \min_{f \in \co f(t, w(t), P)} \partial^\alpha \{\varphi(t, w(\cdot)) \mid f\}
            \leq 0
            \leq \min_{u \in P} \partial^\alpha \{\varphi(t, w(\cdot)) \mid f(t, w(t), u)\}.
        \end{equation*}
    \end{corollary}

    Let us finally assume that the function $f$ from the right-hand side of the dynamic equation \eqref{system} also satisfies the following condition:
    \begin{description}
        \item[$(f.5)$]
            The set $f(t, x, P)$ (see \eqref{f(P)}) is convex for all $t \in [0, T)$ and all $x \in \mathbb{R}^n$.
    \end{description}
    In this case, Corollary \ref{corollary_non-convex} immediately implies the criteria below.
    \begin{corollary} \label{corollary_convex}
        Under assumptions $(f.1)$, $(f.3)$--$(f.5)$, and $(\sigma.2)$, a functional $\varphi \colon G \to \mathbb{R}$ is the value functional of the optimal control problem \eqref{system}, \eqref{cost_functional} if and only if $\varphi$ is continuous, possesses the local Lipschitz continuity property $(L)$, is directionally differentiable of order $\alpha$, meets the boundary condition \eqref{boundary_condition}, and satisfies the equality
        \begin{equation} \label{DI_single-valued_convex}
            \min_{u \in P} \partial^\alpha \{\varphi(t, w(\cdot)) \mid f(t, w(t), u)\}
            = 0
            \quad \forall (t, w(\cdot)) \in G^0.
        \end{equation}
    \end{corollary}

    Let us observe that, according to \eqref{Hamiltonian} and \eqref{directional_derivative_via_ci-derivatives}, at every point $(t, w(\cdot)) \in G^0$ of $ci$-differentiability of order $\alpha$ of a functional $\varphi \colon G \to \mathbb{R}$, the equality from \eqref{DI_single-valued_convex} turns into the equality from \eqref{HJB}.
    In this sense, equation \eqref{DI_single-valued_convex} can be considered as a non-smooth generalization of the Hamilton--Jacobi--Bellman equation \eqref{HJB}.

\section{Optimal Positional Control Strategy in Non-Smooth Case}
\label{section_optimal_control}

    The main result of the present paper is
    \begin{theorem} \label{theorem_optimal_control}
        Let assumptions $(f.1)$, $(f.3)$--$(f.5)$, and $(\sigma.2)$ hold.
        Then, a positional control strategy $U^\circ$ satisfying the condition
        \begin{equation} \label{optimal_strategy_non-smooth}
            U^\circ(t, w(\cdot))
            \in \argmin{u \in P} \partial^\alpha \{ \rho(t, w(\cdot)) \mid f(t, w(t), u) \}
            \quad \forall (t, w(\cdot)) \in G^0
        \end{equation}
        is optimal in the problem \eqref{system}, \eqref{cost_functional}.
        In the above, $\partial^\alpha \{ \rho(t, w(\cdot)) \mid f(t, w(t), u) \}$ is the derivative of order $\alpha$ of the value functional $\rho$ of this problem at the point $(t, w(\cdot))$ in the direction $f(t, w(t), u)$.
    \end{theorem}

    The proof of Theorem \ref{theorem_optimal_control} is based on the properties of the value functional $\rho$ established in Corollary \ref{corollary_convex}, the representation of this functional given by \eqref{value_functional_representation} and \eqref{psi}, formulas \eqref{directional_derivative_formula} and \eqref{N^0} for computing its directional derivatives of order $\alpha$, and the property of uniform directional differentiability of order $\alpha$ of the functional $\psi$ from \eqref{psi} established in Lemma \ref{lemma_psi_is_differentiable}.
    Let us note also that the proof uses the notation introduced throughout the paper (see, in particular, Sections \ref{section_positional} and \ref{subsection_auxiliary_notation}).

    {\it Proof}
        Let us fix $(t, w(\cdot)) \in G^0$ and $\varepsilon > 0$.
        Let us choose $N \geq c_f$ from the condition $(t, w(\cdot)) \in G_N$.
        Using assumption $(f.1)$, let us take $M \geq 0$ such that
        \begin{equation*}
            \|f(\tau, x, u)\|
            \leq M
            \quad \forall \tau \in [0, T] \quad \forall x \in B(R_N) \quad \forall u \in P.
        \end{equation*}
        Due to compactness of the set $X_N$ in $\AC^\alpha([0, T], \mathbb{R}^n)$ and continuity of mapping \eqref{general_mapping} and the value functional $\rho$, there exists $\eta \in (0, T - t)$ such that, for any $\tau \in [T - \eta, T]$ and any $x(\cdot) \in X_N$,
        \begin{equation*}
            |\rho(\tau, x_\tau(\cdot)) - \rho(T, x(\cdot))|
            \leq \frac{\varepsilon}{3}.
        \end{equation*}
        Note that $(t, w(\cdot)) \in G_N^\eta$ by construction.
        Applying Lemma \ref{lemma_psi_is_differentiable}, let us choose $\delta_1 \in (0, \eta / 2]$ such that, for any $(t^\prime, w^\prime(\cdot)) \in G_N^\eta$, any $f \in B(M)$, any $\nu(\cdot) \in \mathcal{N}$, and any $\delta \in (0, \delta_1]$,
        \begin{align*}
            & \biggl| \frac{\psi(t^\prime + \delta, x^{(f)}_{t^\prime + \delta}(\cdot \mid t^\prime, w^\prime(\cdot)), \nu(\cdot))
            - \psi(t^\prime, w^\prime(\cdot), \nu(\cdot))}{\delta} \\
            & \quad - \partial^\alpha_t \psi(t^\prime, w^\prime(\cdot), \nu(\cdot))
            - \langle \nabla^\alpha \psi(t^\prime, w^\prime(\cdot), \nu(\cdot)), f \rangle \biggr|
            \leq \frac{\varepsilon}{3 T}.
        \end{align*}
        Since the value functional $\rho$ possesses property $(L)$, it can be proved similarly to \cite[Proposition 3]{Gomoyunov_Lukoyanov_2021} that there exists $\kappa \geq 0$ such that, for any $(t^\prime, w^\prime(\cdot)) \in G_N^\eta$, any $x(\cdot)$, $x^\prime(\cdot) \in X(t^\prime, w^\prime(\cdot)) \cap X_{N + M}$, and any $\tau \in [t^\prime, T - \eta / 2]$,
        \begin{align*}
            & |\rho(\tau, x_\tau(\cdot)) - \rho(\tau, x^\prime_\tau(\cdot))| \\
            & \quad \leq \kappa \int_{t^\prime}^{\tau} \|(^C D^\alpha x)(\xi) - (^C D^\alpha x^\prime)(\xi)\| \, \rd \xi
            + \kappa (\tau - t^\prime)^{\alpha + 1}.
        \end{align*}
        Using assumption $(f.1)$ and compactness of the set $X_N$ in $\AC^\alpha([0, T], \mathbb{R}^n)$, let us choose $\delta_2 > 0$ such that, for any $\tau$, $\tau^\prime \in [0, T]$ with $|\tau - \tau^\prime| \leq \delta_2$, any $x(\cdot) \in X_N$, and any $u \in P$,
        \begin{equation*}
            \kappa \|f(\tau, x(\tau), u) - f(\tau^\prime, x(\tau^\prime), u)\| + \kappa \delta_2^\alpha
            \leq \frac{\varepsilon}{3 T}.
        \end{equation*}
        Let us put $\delta \coloneqq \min\{\delta_1, \delta_2\} > 0$.

        Let us take a partition $\Delta \coloneqq \{\tau_j\}_{j \in \overline{1, k + 1}}$ of $[t, T]$ satisfying the condition $\diam(\Delta) \leq \delta$ and consider the control $u^\circ(\cdot) \coloneqq u(\cdot \mid t, w(\cdot), U^\circ, \Delta)$ and the motion $x^\circ(\cdot) \coloneqq x(\cdot \mid t, w(\cdot), u^\circ(\cdot))$ of system \eqref{system}, which correspond to the positional control strategy $U^\circ$ from \eqref{optimal_strategy_non-smooth} and the partition $\Delta$.
        Let us note that $x^\circ(\cdot) \in X_N$.
        Since $\tau_1 = t < T - \eta$, let us choose $m \in \overline{2, k + 1}$ such that $\tau_{m - 1} < T - \eta \leq \tau_m$.
        For every $j \in \overline{1, m - 1}$, let us denote $u^\circ_j \coloneqq U^\circ(\tau_j, x^\circ_{\tau_j}(\cdot))$ and $f_j \coloneqq f(\tau_j, x^\circ(\tau_j), u^\circ_j)$.
        We derive
        \begin{align*}
            & \rho(T, x^\circ(\cdot)) - \rho(t, w(\cdot))
            = \rho(T, x^\circ(\cdot)) - \rho(\tau_m, x^\circ_{\tau_m}(\cdot)) \\
            & \quad + \sum_{j = 1}^{m - 1} \bigl( \rho(\tau_{j + 1}, x^\circ_{\tau_{j + 1}}(\cdot))
            - \rho(\tau_{j + 1}, x^{(f_j)}_{\tau_{j + 1}}(\cdot \mid \tau_j, x^\circ_{\tau_j}(\cdot))) \bigr) \\
            & \quad + \sum_{j = 1}^{m - 1} \bigl( \rho(\tau_{j + 1}, x^{(f_j)}_{\tau_{j + 1}}(\cdot \mid \tau_j, x^\circ_{\tau_j}(\cdot)))
            - \rho(\tau_j, x^\circ_{\tau_j}(\cdot)) \bigr),
        \end{align*}
        where the function $x^{(f_j)}(\cdot \mid \tau_j, x^\circ_{\tau_j}(\cdot))$ is defined in accordance with \eqref{x^f}.
        In view of the inclusion $\tau_m \in [T - \eta, T]$, we obtain
        \begin{equation*}
            \rho(T, x^\circ(\cdot)) - \rho(\tau_m, x^\circ_{\tau_m}(\cdot))
            \leq \frac{\varepsilon}{3}.
        \end{equation*}
        Now, let $j \in \overline{1, m - 1}$ be fixed.
        Note that $(\tau_j, x^\circ_{\tau_j}(\cdot)) \in G_N^\eta$ and $f_j \in B(M)$, and, therefore, $x^{(f_j)}(\cdot \mid \tau_j, x^\circ_{\tau_j}(\cdot)) \in X(\tau_j, x^\circ_{\tau_j}(\cdot)) \cap X_{N + M}$.
        Hence, taking into account that $\tau_{j + 1} \leq T - \eta / 2$ and $u^\circ(\tau) = u_j^\circ$ for all $\tau \in [\tau_j, \tau_{j + 1})$, we get
        \begin{align*}
            & \rho(\tau_{j + 1}, x^\circ_{\tau_{j + 1}}(\cdot))
            - \rho(\tau_{j + 1}, x^{(f_j)}_{\tau_{j + 1}}(\cdot \mid \tau_j, x^\circ_{\tau_j}(\cdot))) \\
            & \quad \leq \kappa \int_{\tau_j}^{\tau_{j + 1}} \|f(\xi, x^\circ(\xi), u^\circ_j) - f_j \| \, \rd \xi
            + \kappa (\tau_{j + 1} - \tau_j)^{\alpha + 1} \\
            & \quad \leq \frac{\varepsilon (\tau_{j + 1} - \tau_j)}{3 T}.
        \end{align*}
        Further, by virtue of \eqref{DI_single-valued_convex} (for the value functional $\rho$) and \eqref{optimal_strategy_non-smooth}, we derive
        \begin{equation*}
            \partial^\alpha \{ \rho(\tau_j, x^\circ_{\tau_j}(\cdot)) \mid f_j \}
            = \min_{u \in P} \partial^\alpha \{ \rho(\tau_j, x^\circ_{\tau_j}(\cdot)) \mid f(\tau_j, x^\circ(\tau_j), u) \}
            = 0,
        \end{equation*}
        and, consequently, due to \eqref{directional_derivative_formula}, there exists $\nu_j^\circ(\cdot) \in \mathcal{N}^\circ(\tau_j, x^\circ_{\tau_j}(\cdot))$ such that
        \begin{equation*}
            \partial^\alpha_t \psi(\tau_j, x^\circ_{\tau_j}(\cdot), \nu_j^\circ(\cdot))
            + \langle \nabla^\alpha \psi(\tau_j, x^\circ_{\tau_j}(\cdot), \nu_j^\circ(\cdot)), f_j \rangle
            = 0.
        \end{equation*}
        Then, in accordance with \eqref{value_functional_representation}, \eqref{psi}, and \eqref{N^0}, we obtain
        \begin{align*}
            & \rho(\tau_{j + 1}, x^{(f_j)}_{\tau_{j + 1}}(\cdot \mid \tau_j, x^\circ_{\tau_j}(\cdot)))
            - \rho(\tau_j, x^\circ_{\tau_j}(\cdot)) \\
            & \quad \leq \psi(\tau_{j + 1}, x^{(f_j)}_{\tau_{j + 1}}(\cdot \mid \tau_j, x^\circ_{\tau_j}(\cdot)), \nu_j^\circ(\cdot))
            - \psi(\tau_j, x^\circ_{\tau_j}(\cdot), \nu_j^\circ(\cdot)) \\
            & \quad \leq \frac{\varepsilon (\tau_{j + 1} - \tau_j)}{3 T}.
        \end{align*}
        As a result, recalling that the value functional $\rho$ satisfies the boundary condition \eqref{value_functional_T}, we conclude that
        \begin{equation*}
            \sigma(x^\circ(T))
            = \rho(T, x^\circ(\cdot))
            \leq \rho(t, w(\cdot)) + \varepsilon,
        \end{equation*}
        which implies that the control $u^\circ(\cdot)$ is $\varepsilon$-optimal and completes the proof.
    \qed

    Finally, let us observe that, according to \eqref{directional_derivative_via_ci-derivatives}, at every point $(t, w(\cdot)) \in G^0$ of $ci$-differentiability of order $\alpha$ of the value functional $\rho$, the inclusion from \eqref{optimal_strategy_non-smooth} turns into the inclusion from \eqref{U^0_smooth}.
    Hence, the rule for constructing an optimal positional control strategy \eqref{optimal_strategy_non-smooth} can be considered as a non-smooth generalization of the rule \eqref{U^0_smooth}.

\section{Example}
\label{section_Example}

    Let us illustrate the results presented in this paper by a model example.
    Let us consider an optimal control problem for a dynamical system described by the Caputo fractional differential equation (we take $n = 1$)
    \begin{equation} \label{example_system}
        (^C D^\alpha x)(\tau)
        = \mathrm{\Gamma}(\alpha) g(\tau) u(\tau),
    \end{equation}
    where $\tau \in [0, T]$, $x(\tau) \in \mathbb{R}$, $u(\tau) \in [- 1, 1]$, and $g \colon [0, T] \to \mathbb{R}$ is a given continuously differentiable function (the multiplier $\mathrm{\Gamma}(\alpha)$ is added for convenience), and the terminal cost functional
    \begin{equation} \label{example_cost_functional}
        J(t, w(\cdot), u(\cdot))
        \coloneqq - \bigl( x(T \mid t, w(\cdot), u(\cdot)) \bigr)^2
    \end{equation}
    for all $(t, w(\cdot)) \in G^0$ and all $u(\cdot) \in \mathcal{U}(t, T)$.
    In this problem, assumptions $(f.1)$, $(f.3)$--$(f.5)$, and $(\sigma.2)$ are clearly hold.
    Let us apply Theorem \ref{theorem_optimal_control} in order to construct an optimal positional control strategy.

    Let us consider the functional
    \begin{equation*}
        \varphi(t, w(\cdot))
        \coloneqq - \biggl( |a(T \mid t, w(\cdot))| + \int_{t}^{T} \frac{|g(\tau)|}{(T - \tau)^{1 - \alpha}} \, \rd \tau \biggr)^2
        \quad \forall (t, w(\cdot)) \in G,
    \end{equation*}
    where the function $a(\cdot \mid t, w(\cdot))$ is given by \eqref{a}.
    Using Corollary \ref{corollary_convex}, let us prove that $\varphi$ is actually the value functional of the problem \eqref{example_system}, \eqref{example_cost_functional}.
    According to \cite[Section 12]{Gomoyunov_2020_SIAM}, the auxiliary functional $\varphi_\ast(t, w(\cdot)) \coloneqq a(T \mid t, w(\cdot))$ for all $(t, w(\cdot)) \in G$ is $ci$-smooth of order $\alpha$ and
    \begin{equation*}
        \partial_t^\alpha \varphi_\ast(t, w(\cdot))
        = 0,
        \quad \nabla^\alpha \varphi_\ast(t, w(\cdot))
        = \frac{1}{\mathrm{\Gamma}(\alpha) (T - t)^{1 - \alpha}}
        \quad \forall (t, w(\cdot)) \in G^0.
    \end{equation*}
    Based on these properties, it can be shown that the functional $\varphi$ is continuous and possesses property $(L)$.
    Moreover, arguing similarly to \cite[Section 12]{Gomoyunov_2020_SIAM}, we derive that, at every point $(t, w(\cdot)) \in G^0$ such that $a(T \mid t, w(\cdot)) \neq 0$, the functional $\varphi$ is $ci$-differentiable of order $\alpha$ and
    \begin{align}
        \partial^\alpha_t \varphi(t, w(\cdot))
        & = \frac{2 |g(t)|}{(T - t)^{1 - \alpha}}
        \biggl( |a(T \mid t, w(\cdot))|
        + \int_{t}^{T} \frac{|g(\tau)|}{(T - \tau)^{1 - \alpha}} \, \rd \tau \biggr),
        \nonumber \\
        \nabla^\alpha \varphi(t, w(\cdot))
        & = \frac{- 2 \sgn (a(T \mid t, w(\cdot)))}{\mathrm{\Gamma}(\alpha) (T - t)^{1 - \alpha}}
        \biggl( |a(T \mid t, w(\cdot))|
        + \int_{t}^{T} \frac{|g(\tau)|}{(T - \tau)^{1 - \alpha}} \, \rd \tau \biggr).
        \label{example_ci_derivatives}
    \end{align}
    In particular, at all such points, the Hamilton--Jacobi--Bellman equation corresponding to the optimal control problem \eqref{example_system}, \eqref{example_cost_functional} is satisfied:
    \begin{equation*}
        \partial^\alpha_t \varphi(t, w(\cdot))
        - \mathrm{\Gamma(\alpha)} |g(t)| |\nabla^\alpha \varphi(t, w(\cdot))|
        = 0.
    \end{equation*}
    At the same time, for every point $(t, w(\cdot)) \in G^0$ with $a(T \mid t, w(\cdot)) = 0$, we get
    \begin{equation} \label{example_directional_derivatives}
        \partial^\alpha \{ \varphi(t, w(\cdot)) \mid f \}
        = \frac{2 (\mathrm{\Gamma}(\alpha) |g(t)| - |f|)}{\mathrm{\Gamma}(\alpha) (T - t)^{1 - \alpha}}
        \int_{t}^{T} \frac{|g(\tau)|}{(T - \tau)^{1 - \alpha}} \, \rd \tau
        \quad \forall f \in \mathbb{R},
    \end{equation}
    and, in particular, the equality below is valid:
    \begin{equation*}
        \min_{u \in [- 1, 1]} \partial^\alpha \{ \varphi(t, w(\cdot)) \mid \mathrm{\Gamma}(\alpha) g(t) u \}
        = 0.
    \end{equation*}
    Thus, taking into account that
    \begin{equation*}
        \varphi(T, w(\cdot))
        = - \bigl( a( T \mid t, w(\cdot)) \bigr)^2
        = - \bigl( w(T) \bigr)^2
        \quad \forall w(\cdot) \in \AC^\alpha([0, T], \mathbb{R}),
    \end{equation*}
    we conclude that $\varphi$ is the value functional of the problem \eqref{example_system}, \eqref{example_cost_functional}.

    Let us emphasize that, according to \eqref{directional_derivative_via_ci-derivatives}, it follows from \eqref{example_directional_derivatives} that, at every point $(t, w(\cdot)) \in G^0$ such that $a(T \mid t, w(\cdot)) = 0$, the functional $\varphi$ is not $ci$-differentiable of order $\alpha$ if the function $g$ is not identically zero on $[t, T]$.

    As a result, using formulas \eqref{example_ci_derivatives} and \eqref{example_directional_derivatives} for the corresponding derivatives of the value functional of the problem \eqref{example_system}, \eqref{example_cost_functional}, we obtain by Theorem \ref{theorem_optimal_control} that the following positional control strategy $U^\circ$ is optimal:
    \begin{equation*}
        U^\circ(t, w(\cdot))
        \coloneqq \begin{cases}
            \sgn (g(t)),
            & \mbox{if } a(T \mid t, w(\cdot)) > 0, \\
            1,
            & \mbox{if } a(T \mid t, w(\cdot)) = 0, \\
            - \sgn (g(t)),
            & \mbox{if } a(T \mid t, w(\cdot)) < 0,
        \end{cases}
        \quad \forall (t, w(\cdot)) \in G^0.
    \end{equation*}

\begin{acknowledgements}
    This work was supported by RSF, project no. 19-11-00105.
\end{acknowledgements}

\end{document}